\documentclass[11pt]{article}

\usepackage{amsmath,amssymb,amsfonts,enumitem,amsthm,accents,pstricks,pst-node,pstricks-add,mathrsfs,xy}
\xyoption{all}

\usepackage{color}
\usepackage{makeidx}
\usepackage{graphicx}
\usepackage{setspace}
\usepackage{extarrows}

\numberwithin{equation}{section}

\usepackage{tocloft}

\newtheorem{theorem}{Theorem}[section]
\newtheorem{lemma}[theorem]{Lemma}
\newtheorem{conjlemma}[theorem]{Conjectural Lemma}
\newtheorem{corollary}[theorem]{Corollary}
\newtheorem{proposition}[theorem]{Proposition}
\newtheorem{conjecture}[theorem]{Conjecture}

\theoremstyle{definition}
\newtheorem{definition}[theorem]{Definition}
\newtheorem{remark}[theorem]{Remark}
\newtheorem{example}[theorem]{Example}

\newtheorem{claim}[theorem]{Claim}

\def\ov#1{\overline{#1}}
\def\tn#1{\textnormal{#1}}
\def\mf#1{\mathfrak{#1}}
\def\wt#1{\widetilde{#1}}

\def\scz{\scriptsize}
\def\vr{\varrho}

\def\mc{\mathcal}
\def\lra{\longrightarrow}
\def\dbar{\bar\partial}

\def\ve{\varepsilon}

\newcommand\uvec[1]{\underaccent{\vec}{#1}}
\newcommand{\ucev}[1]{\reflectbox{\ensuremath{\uvec{\reflectbox{\ensuremath{#1}}}}}}

\newcommand{\abs}[1]{\left\vert #1 \right\vert}

\newcommand{\lrc}[1]{\left\{ #1 \right\}}

\newcommand*{\defeq}{\mathrel{\vcenter{\baselineskip0.5ex \lineskiplimit0pt
                     \hbox{\scriptsize.}\hbox{\scriptsize.}}}%
                     =}

\def\scz{\scriptsize}

\def\bEq#1{\begin{equation}\label{#1}}
\def\eEq{\end{equation}}

\def\bsEq{\begin{equation*}}
\def\esEq{\end{equation*}}

\def\bDf#1{\begin{definition}\label{#1}}
\def\eDf{\end{definition}}

\def\bTh#1{\begin{theorem}\label{#1}}
\def\eTh{\end{theorem}}

\def\bCn#1{\begin{conjecture}\label{#1}}
\def\eCn{\end{conjecture}}

\def\bLm#1{\begin{lemma}\label{#1}}
\def\eLm{\end{lemma}}

\def\bCLm#1{\begin{conjlemma}\label{#1}}
\def\eCLm{\end{conjlemma}}

\def\bRm#1{\begin{remark}\label{#1}}
\def\eRm{\end{remark}}

\def\bEx#1{\begin{example}\label{#1}}
\def\eEx{\end{example}}

\def\bPr#1{\begin{proposition}\label{#1}}
\def\ePr{\end{proposition}}

\def\bCr#1{\begin{corollary}\label{#1}}
\def\eCr{\end{corollary}}

\def\bFg#1{\begin{figure}\label{#1}}
\def\eFg{\end{figure}}

\def\bPf{\begin{proof}}
\def\ePf{\end{proof}}

\def\bIt{\begin{itemize}[leftmargin=*]}
\def\eIt{\end{itemize}}

\def\bEn{\begin{enumerate}[label=$(\arabic*)$,leftmargin=*]}
\def\eEn{\end{enumerate}}

\def\bEnalph{\begin{enumerate}[label=$(\alph*)$,leftmargin=*]}
\def\eEnalph{\end{enumerate}}

\def\nd{\tn{d}}

\def\Def{\tn{Def}}

\def\Obs{\tn{Obs}}
\def\ord{\tn{ord}}

\def\sing{\tn{sing}}

\def\vfc{\tn{VFC}}

\def\ob{\tn{ob}}

\def\cM{\mc{M}}
\def\cO{\mc{O}}

\def\cP{\mc{P}}

\def\cN{\mc{N}}

\def\cY{\mc{Y}}
\def\cZ{\mc{Z}}

\def\cL{\mc{L}}

\def\cT{\mc{T}}
\def\cI{\mc{I}}
\def\E{\mathbb E}
\def\R{\mathbb R}
\def\C{\mathbb C}
\def\Z{\mathbb Z}

\def\P{\mathbb P}

\def\T{\mathbb T}
\def\N{\mathbb N}
\def\D{\mathbb D}
\def\L{\mathbb L}

\def\V{\mathbb V}
\def\mfi{\mf{i}}
\def\mfj{\mf{j}}

\def\aut{\mf{Aut}}

\def\mfs{\mf{s}}

\def\la{\lambda}
\def\La{\Lambda}

\def\De{\Delta}
\def\de{\delta}

\def\Om{\Omega}
\def\si{\sigma}
\def\Si{\Sigma}
\def\al{\alpha}

\def\Ga{\Gamma}
\def\ze{\zeta}

\def\eset{\emptyset}

\makeindex
\begin{document}
\title{Limits of stable maps in a semi-stable degeneration}
\author{Mohammad Farajzadeh Tehrani\footnote{This work is supported by the NSF grant DMS-2003340} }
\date{\today}
\maketitle

\begin{abstract}
Given a semistable degeneration with a  simple normal crossings central fiber, Abramovich-Chen-Gross-Siebert \cite{ACGS} proved a degeneration formula that relates the moduli spaces of stable maps  in smooth fibers to certain moduli spaces of log-smooth maps in the central fiber. 
In this paper, we study the same problem from an analytic point of view. We prove that  the limiting stable maps in the central fiber satisfy specific combinatorial and analytical conditions.
Furthermore, we explain the deformation-obstruction theory of the moduli spaces arising from these conditions, derive a degeneration formula, and work out an explicit example.
The earlier version \cite{FTowards} of this paper contains an outline of these ideas for the symplectic category.\end{abstract}
\tableofcontents

%------------------------------------------------------------------------------------------------------
\section{Introduction}\label{intro_s}
In this paper, by a \textbf{semistable degeneration} we mean a one-parameter family 
\bEq{pi_e}
\pi\!:\cZ\!\lra\!\De,
\eEq
where $\De$ is a  compact disk around the origin in $\C$, $\cZ$ is a smooth K\"ahler manifold with compact fibers, the central fiber 
$$
\cZ_0\!\defeq\!\pi^{-1}(0)=X_\eset\defeq\bigcup_{i\in \cI} X_i
$$ 
is a compact simple normal crossings (or \textbf{SNC}) K\"ahler variety, and the fibers over $\De^*\!\defeq\!\De\!-\!\{0\}$ are smooth. An SNC variety with 3 irreducible components is shown in Figure~\ref{3conf_fig}. 
\begin{figure}
\begin{pspicture}(-3,-2)(11,.5)
\psset{unit=.3cm}
\psline[linewidth=.1](15,-2)(21,-2)\psline[linewidth=.1](15,-2)(15,4)
\psline[linewidth=.1](15,-2)(11.34,-5.66)\pscircle*(15,-2){.3}
\rput(19.5,2.5){\small{$X_i$}}\rput(11.5,-1){\small{$X_j$}}\rput(17,-5){\small{$X_k$}}
\rput(22.2,-2.1){\small{$X_{ik}$}}\rput(15.1,4.8){\small{$X_{ij}$}}
\rput(10.8,-6.1){\small{$X_{jk}$}}\rput(16.3,-1.2){\small{$X_{ijk}$}}
\end{pspicture}
\caption{A 3-fold SC variety.}
\label{3conf_fig}
\end{figure}
For each $i\!\in\!\cI$, let $\cN_i\!\equiv\!\cN_\cZ X_i$ denote the normal line bundle of $X_i$ in $\cZ$.
The line bundle $\cO_{\cZ}(\cZ_0)$ is trivial. Let
$$
X_I\!\equiv\!\bigcap_{i\in I} X_i\qquad \forall~\eset\neq I\subset \mc{I}.
$$
Any trivialization $\cO_{\cZ}(\cZ_0)\cong \cO_{\cZ}$ restricts to a  set of compatible trivializations 
\bEq{LBC_e}
\cO_{\cZ}(\cZ_0)|_{X_I}=\bigotimes_{i\in \mc{I}}\cO_{\cZ}(X_i)|_{X_{I}}=
\bigotimes_{i\in I}\cN_i|_{X_{I}} \otimes \bigotimes_{i\in \cI-I} \cO_{X_I}(X_{I+i})\cong \cO_{X_I},
\eEq
where $\cO_{X_I}(X_{I+ i})$ is the line bundle corresponding to the smooth divisor 
$$
X_{I+i}\defeq X_{I\cup \{i\}}\!\subset\!X_I, \qquad \forall~i\!\in\!\cI-I.
$$
The trivialization (\ref{LBC_e}) can also be obtained from the $|I|$-th derivative of the projection map $\pi$ in (\ref{pi_e}) along $X_I$.\\

\noindent
For every $\la\!\in\!\De^*$, given $g,k\!\in\!\N$ and  $A\!\in\!H_2(\cZ_\la,\Z)$, a $k$-marked genus $g$ degree $A$ stable map into $\cZ_\la$ is a tuple 
$$
\big(u\colon\!\Si\!\lra\!\cZ_\la, \vec{z}=(z_1,\ldots,z_k)\big)
$$  
where $C\!\equiv\!(\Si,\vec{z})$ is a connected genus $g$ nodal curve with $k$ distinct ordered marked points (away from the nodes) and $u$ is a holomorphic map representing the homology class $A$. 
Two marked stable maps 
$$
\big(u\colon\!\Si\!\lra\!\cZ_\la, \vec{z}\big)\quad\tn{and}\quad\big(\wt{u}\colon\!\wt{\Si}\!\lra\!\cZ_\la, \vec{\wt{z}}\big)
$$
are \textbf{equivalent} if there exists a bi-holomorphic isomorphism $h\colon\!\Si\!\lra\!\wt{\Si}$ such that $h(z_{a})\!=\!\wt{z}_{a}$, for all $a\!=\!1,\ldots,k$, and $u\!=\!\wt{u}\circ h$. A marked stable map is \textbf{stable} iff the group of self-automorphisms is finite. For $\la\!\in\!\De^*$, let $\ov\cM_{g,k}(\cZ_\la,A)$ denote the moduli space (set) of equivalence classes of $k$-marked genus $g$ degree $A$ stable maps into $\cZ_\la$.\\

\noindent
If $\tn{dim}_\C\cZ_\la=n$, the expected $\C$-dimension of $\ov\cM_{g,k}(\cZ_\la,A)$ is
\bEq{exp-dim_e}
c_1^{T\cZ_\la}(A)+(n-3)(1-g)+k.
\eEq
Gromov-Witten (or \textbf{GW}) invariants are obtained by the integration of certain cohomology classes against the virtual fundamental class (or \textbf{VFC}) of $\ov\cM_{g,k}(\cZ_\la,A)$. Since there might be different homology classes in $\cZ_\la$ that are the same as homology classes in $\cZ$, for each $A\!\in\!H_2(\cZ,\Z)$, we let $\ov\cM_{g,k}(\cZ_\la,A)$ to be the union over all the representatives of $A$ in $H_2(\cZ_\la,\Z)$; see \cite{FZRim} for a careful discussion of this issue. \\

\noindent
For any choice of $(g,k,A)$, the fiberation (\ref{pi_e}) gives rise to a $1$-parameter family  
\bEq{Family*_e}
\ov\cM_{g,k}(\cZ^*,A)\defeq \bigcup_{\la\in\De^*}\ov\cM_{g,k}(\cZ_\la,A)\lra \De^*
\eEq
with fibers of equal virtual dimension (and cobordant VFC). For $\la\!=\!0$, let $\ov\cM_{g,k}(\cZ_0,A)$ denote the space of all stable maps in $\cZ$ whose image lies inside $\cZ_0$; $\ov\cM_{g,k}(\cZ_0,A)$ is not a moduli space of the correct expected dimension that extends the virtual cobordism (\ref{Family*_e}) over $0\in \De$. Therefore, from an analytical perspective, the important questions are:
\bEn
\item\label{Q1}
\textit{which stable maps in $\ov\cM_{g,k}(\cZ_0,A)$ can (virtually) arise as the Gromov-limit of a sequence of stable maps in~(\ref{Family*_e})? }
\item\label{Q2}
\textit{how to complete (\ref{Family*_e}) with a moduli space $\ov\cM^{\tn{good}}_{g,k}(\cZ_0,A)$, ideally still a subset of $\ov\cM_{g,k}(\cZ_0,A)$, admitting a VFC that is cobordant to VFC of smooth fibers?} 
\item\label{Q3}
\textit{can $\ov\cM^{\tn{good}}_{g,k}(\cZ_0,A)$ and its VFC be expressed in terms of certain moduli spaces in $X_i$ relative to the SNC divisor $\partial X_i\equiv \bigcup_{j\in \cI-i} X_{ij} \subset X_i$ and their VFCs?} 
\eEn

\noindent
In the algebraic category and for a semistable degeneration into two pieces $\cZ_0\!=\!X_1\!\cup_{X_{12}}\!X_2$ along a smooth divisor,  these questions were first answered by Jun Li \cite{JLi1,JLi2}. For a smooth divisor $D\!\subset\!X$, he introduced the notion of a \textit{stable relative map} whose image lives in a natural SNC ``expanded degeneration" associated to $(X,D)$. Similarly, for a semistable degeneration  into two pieces $\cZ_0\!=\!X_1\!\cup_{X_{12}}\!X_2$,  he constructed a compactification $\ov\cM^{\tn{rel}}_{g,k}(\cZ_0,A)$ whose (virtually) main components are fiber products of the relative moduli spaces\footnote{Over possibly disconnected domains with Euler characteristic $\chi_1$ and $\chi_2$.} $\ov\cM_{\chi_1,\mfs}(X_1/X_{12},A_1)$ and $\ov\cM_{\chi_2,\mfs}(X_2/X_{12},A_2)$. In \cite{JLi2}, he proved a decomposition formula which expresses the GW invariants of the smooth fibers in terms of the products of relative GW invariants of $(X_1,X_{12})$ and $(X_2,X_{12})$. For a symplectic version of these results see \cite{LR,IPsum,FZSum}.
More recently, Gross-Siebert \cite{GS} and Abramovich-Chen \cite{ACII} introduced moduli spaces of (fine, saturated) stable log maps and proved a degeneration formula \cite{ACGS} to answer the first two questions above in an arbitrary semistable degeneration. Also, in \cite{PFormula}, Brett Parker uses moduli spaces of curves in exploded manifolds \cite{BP1,BP2, BP3} to address the first two questions. These constructions work for even a larger class of ``log smooth" (see \cite{ACGHOSS}) and ``exploded" varieties (see \cite{P-Exp}), respectively. The degeneration formula\footnote{or as they call it: the ``invariance property".} in \cite{ACGS} can be read as: \textit{virtually, a stable map $f$ in $\ov\cM_{g,k}(\cZ_0,A)$ can arise as a limit of a sequence of  stable maps in smooth fibers if and only if  $f$ can be enhanced to a log smooth map in the log moduli space\footnote{The notation $\ov\cM^{\tn{al,}\log}_{g,k}(\cZ_0,A)$ for the algebraic log moduli space is not the one in \cite{ACII,GS}. } $\ov\cM^{\tn{al},\log}_{g,k}(\cZ_0,A)$}. A log smooth map is a stable map plus a lift of that to a map between certain sheaves of monoids satisfying some conditions. Parker's definition in the category of exploded manifolds involves sheaf theory in a similar way. For a  geometric approach using expanded degenerations we refer to \cite{R}. For the explanation of the difference between Jun Li's formula and ACGS's formula we refer to \cite{KLR}. \\

\noindent
Theorem~5.3.3 in \cite{ACGS}, gives a criterion for lift-ability and a formula for the number of lifts.
However, in practice, given a stable map $f$, its is rather hard to check whether $f$ lifts to an element of $\ov\cM^{\tn{al},\log}_{g,k}(\cZ_0,A)$. In other words, the image of the forgetful map 
$$
\ov\cM^{\tn{al},\log}_{g,k}(\cZ_0,A)\lra \ov\cM_{g,k}(\cZ_0,A)
$$
is hard to describe. For example, in \cite[Sec.~6.2.5]{ACGS}, the authors lift to a blowup of the SNC central fiber to find those ``star-shaped" maps that can be lifted to a log map.  \\

\noindent
The goal of this paper is two-fold. First, we prove that the Gromov-limits of stable maps in the central fiber satisfy two specific combinatorial (called C1) and analytical (called C2) conditions. The combinatorial condition (C1) is equivalent to the basicness  condition in \cite[Dfn.~1.20]{GS}. The analytical condition (C2) has no direct analogue in \cite{ACGS,PFormula}. The two conditions are linked by a linear map associated to the dual graph of the stable map in the question. The positive cone in the kernel of this map gives a toric description of the space of gluing parameters. It is interesting and important to figure out the relation between (C2) and the lift-ability criterion in \cite[Thm.~5.3.3]{ACGS}. We plan to work on this in the future. Second, we provide evidence that the moduli space of log maps satisfying conditions (C1) and (C2) should similarly address the first two questions above. In particular, we derive an explicit degeneration formula that, in the case of basic degenerations, coincides with Jun Li's formula.
The degeneration formula \cite[(1.1.1)]{ACGS} and the one that we propose here are both a sum over the same set of combinatorial data, but with different coefficients; see Remark~\ref{ACG_rmk}. The only obstacle in the way of generalizing these result to the symplectic category is to find a suitable class of almost complex structures compatible with a symplectic semistable degeneration (in the sense of \cite{FMZSum}); see \cite{FTowards}. 

\begin{remark}
Since $\cZ_0\!\subset\! \cZ$ is an SNC divisor, Theorem 1.3 in \cite{FRelative} with trivial tangency data at the marked points gives us a compact (relative) log moduli space $\ov\cM^{\log}_{g,k}(\cZ,\cZ_0,A)$ that contains $\ov\cM_{g,k}(\cZ^*,A)$ as an open subset. However, if $g\!>\!0$, even the expected dimension of the subset of log curves in $\ov\cM^{\log}_{g,k}(\cZ,\cZ_0,A)$ that live in $\cZ_0$ is different from (\ref{exp-dim_e}). The conditions (C1) and (C2) are a refinement of the similar conditions in \cite[Dfn.~2.8]{FRelative}. Therefore, the compactness result \cite[Thm.~1.3]{FRelative} does not directly apply; it needs some enhancements.
\end{remark}

\noindent
For a finite set $\cI$ and a ring $R$, let
$$
R^\cI_\bullet=\big\{r=(r_j)_{j\in \cI}\colon \sum_{j\in \cI} r_j=0\big\}\subset R^\cI.
$$ 
For each $i\!\in\!\cI$, let $\xi_{i}$ denote a non-zero holomorphic section of $\cO_{\cZ}(X_i)$ vanishing (to the order $1$) along $X_i$. Since $\cZ_0$ is compact, the restriction of each section $\xi_i$ to $\cZ_0$ is unique up to multiplication by a constant. We will choose these sections so that the composition 
\bEq{xiChoice_e}
\cZ \lra \cO_\cZ(\cZ_0) \cong \cZ\times \C \lra \C,
\eEq
where the first map is $x\to \prod_{i\in \mc{I}} \xi_i(x)$ and the last map is projection to the second factor, is equal to (\ref{pi_e}).
For the trivial holomorphic line bundle $\cO$ (on any base), let $\la_{\cO}$ denote the constant section corresponding to $\la\!\in\! \C$.\\ 

\noindent
Fix a trivialization of $\cO_{\cZ}(\cZ_0)$. With notation as above, we define an analytical marked nodal log map into $\cZ_0$ with the marked nodal domain $(\Si,\vec{z})\!=\!\bigcup_{v\in \V}(\Si_v,\vec{z}_v)$ to be a collection of tuples
\bEq{logmap_first}
f\equiv \big(u_v\colon\!\Si_v\!\lra\! X_{I_v}, \vec{z}_v, (\ze_{v,i})_{i\in I_v}\big)_{v\in \V}
\eEq
over smooth components $\Si_v$ of $\Si$ such that 
\bEn
\item $\big(u\!\equiv\!(u_v)_{v\in \V}\colon\! \Si\!\lra\! \cZ_0,\vec{z}\big)$ is a $k$-marked nodal map in the classical sense,
\item for each $v\!\in\!\V$, $\eset\!\neq\!I_v\!\subset\!\cI$ is the maximal subset such that $\tn{Im}(u_v)\!\subset\!X_{I_v}\!\subset\!\cZ_0$,
\item for each $v\!\in\!\V$ and every $i\!\in\!I_v$, $\ze_{v,i}$ is a non-trivial meromorphic section of the holomorphic line bundle $u_v^* \cN_i$,
\item for each $v\in\V$, with respect to the isomorphism (\ref{LBC_e}),  we have
$$
\bigotimes_{i\in I} \ze_{v,i} \otimes \bigotimes_{j\in \cI-I} u_v^* \xi_{j} = u_v^*(\la_{\cO_{\cZ}(\cZ_0)})
$$
for some fixed $\la\!=\!\la(f)\!\in\!\C^*$ (independent of $v\in \V$),
\item the ``contact order vectors" in $\Z_\bullet^\cI$, defined in (\ref{Ordx_e1}) and  (\ref{Ordx_e2}), are the opposite of each other at the nodal points of $\Si$, 
\item every point in $\Si$ with a non-trivial contact vector is a nodal point,
\item (C1:) there exists a vector-valued function $s\colon\!\V\!\lra\!\R^\cI$ such that $s_v\!=\!s(v)\!\in\! \R_{+}^{I_v}\!\times \{0\}^{\cI-I_v}$ for all $v\!\in\!\V$, and $s_v\!-\!s_{v'}$ is a positive multiple of the contact order vector of any nodal point on $\Si_v$ connected to $\Si_{v'}$, for all $v,v'\!\in\!\V$, 
\item\label{GObs_it}  (C2:) certain Lie group (a complex torus) element $\ob(f)$ associated to $f$, defined in (\ref{PLtoG_e}), is equal to $1$;
\eEn
see Definitions~\ref{LogTuple_dfn}, \ref{PreLogMap_dfn}, and  \ref{LogMapVar_dfn} for the details. In simple words, a log map is a stable map together with a set of meromorphic sections that satisfies certain combinatorial (i.e. (5), (6), and (C1)) and analytical (i.e. (4) and (C2)) conditions.\\

\noindent
Two marked log maps 
$$
f\equiv \big(u_v\colon\!\Si_v\!\lra\! X_{I_v}, \vec{z}_v, (\ze_{v,i})_{i\in I_v}\big)_{v\in \V}\quad \tn{and}\quad
\wt{f}\equiv \big(\wt{u}_v\colon\!\wt{\Si}_v\!\lra\! X_{I_v}, \vec{\wt{z}}_v, (\wt{\ze}_{v,i})_{i\in I_v}\big)_{v\in \V}
$$
are \textbf{equivalent} if there exists a bi-holomorphic isomorphism
$$
(h\colon \Si\lra \wt{\Si})\equiv \big(h_v\colon \Si_v\lra\wt\Si_{h(v)}\big)_{v\in \V}
$$ 
such that 
$$
h(z_{a})=\wt{z}_{a}\quad \forall~a=1,\ldots,k, \qquad \wt{u}\circ h= u, \qquad  h_v^*\wt{\ze}_{h(v),i}=c_{v,i} \ze_{v,i} \quad \forall~v\!\in\!\V,~i\!\in\!I_v.
$$
In particular, given a marked log map $f$ as in (\ref{logmap_first}), replacing each meromorphic section $\ze_{v,i}$ with a non-zero multiple $c_{v,i} \ze_{v,i}$ of that, satisfying $\prod_{i\in I_v}{c_{v,i}}=c'$ for all $v\!\in\! \V$, produces another marked log map which is equivalent to $f$. 
A marked log map is stable if it has a finite automorphism group. For $g,k\!\in\!\N$ and $A\!\in\!H_2(\cZ_0,\Z)$, \textit{we denote the space of equivalence classes of stable  $k$-marked degree $A$ genus $g$  log maps by} 
$$
\ov\cM_{g,k}^{\log}(\cZ_0,A). 
$$
This moduli space is independent of the choice of the sections $\xi_i$ used in the construction because rescalings of $\xi_i$ can be compensated by rescalings of $\ze_{v,i}$. 
The equivalence class of an analytic log map is called an analytic log \textit{curve}. We will often drop the adjective ``analytic" and simply say log map or log curve.\\ 

\noindent
There is a natural forgetful map
$$
\ov\cM_{g,k}^{\log}(\cZ_0,A)\lra \ov\cM_{g,k}(\cZ_0,A),\qquad \big(u_v\colon\!\Si_v\!\to\! X_{I_v},\vec{z}_v, (\ze_{v,i})_{i\in I_v}\big)_{v\in \V}\lra \big(u_v\colon\!\Si_v\!\to\! \cZ_0,\vec{z}_v\big)_{v\in \V}.
$$
It turns out  that for every $k$-marked stable nodal curve $f$ in $\ov\cM_{g,k}(\cZ_0,A)$, there exists at most finitely many log curves $f_{\log}\!\in\!\ov\cM^{\log}_{g,k}(\cZ_0,A)$ (with distinct decorations on the dual graph) lifting $f$; see Remark~\ref{elboration_rmk}. Furthermore, $f_{\log}$ is stable if and only if $f$ is stable, and the automorphism groups are often the same. 

\bRm{nocZ_rmk}
In Section~\ref{LogSCVar_s}, we will construct the analytical log moduli spaces for any arbitrary $d$-semistable (see \cite{Fr})  SNC variety $\cZ_0$ without using the smoothing $\cZ$ that contains it. Here, we used $\cZ$ to slightly simplify the notation. Furthermore, it is possible to define the log map without mentioning the meromorphic sections $\ze_{v,i}$; see Remark~\ref{elboration_rmk}.
\eRm

\begin{definition}
A continuous function $f\colon\! M\!\lra\! N$ between two topological spaces is a local embedding if for all $x\!\in\!M$ there is an open neighborhood $U\!\ni\!x$ such that $f|_U\colon\! U\!\lra\! N$ is an embedding.
\end{definition}

\noindent
By Smirnov's theorem, every paracompact, Hausdorff, and locally metrizable space is metrizable.
Therefore, if $f\colon\! M\!\lra\! N$ is a local embedding from a compact Hausdorff space $M$ to a compact metrizable space $N$ then $M$ is metrizable.

\bTh{Compactness_th}
For every $A\!\in\!H_2(\cZ_0,\Z)$ and $g,k\!\in\!\N$, the Gromov sequential convergence topology on $\ov\cM_{g,k}(\cZ,A)$ lifts to a compact Hausdorff sequential convergence topology
on 
$$
 \ov\cM^{\log}_{g,k}(\cZ,A)\defeq \ov\cM_{g,k}(\cZ^*,A) \cup \ov\cM^{\log}_{g,k}(\cZ_0,A)
$$ 
such that the natural forgetful maps
\bEq{FogetLog_e}
 \iota\colon \ov\cM^{\log}_{g,k}(\cZ,A)\lra \ov\cM_{g,k}(\cZ,A)
\quad \tn{and}\quad 
\iota \colon \ov\cM^{\log}_{g,k}(\cZ_0,A)\lra \ov\cM_{g,k}(\cZ_0,A)
\eEq
are local embeddings. In particular, $\ov\cM^{\log}_{g,k}(\cZ,A)$ and $\ov\cM^{\log}_{g,k}(\cZ_0,A)$ are metrizable.
If $g\!=\!0$, then the forgetful maps in (\ref{FogetLog_e}) are global embeddings.
\eTh

\noindent
If $X_\eset$ is just an abstract d-semistable SNC variety, we just get the restriction of the theorem above to $\ov\cM^{\log}_{g,k}(X_\eset,A)$.
If $\cZ_0$ is basic ($\cZ_0\!=\!X_1\cup_{X_{12}} X_2$), it follows from \cite[Prp.~4.5]{FRelative} that there is a surjective projection map 
$$
\ov\cM^{\tn{rel}}_{g,k}(\cZ_0,A)\!\lra\!\ov\cM^{\log}_{g,k}(\cZ_0,A),
$$ 
where the former is Jun Li's relative moduli space. The degeneration formula that we will drive will be the same as Jun Li's formula in this case.\\

\noindent
Theorem~\ref{Compactness_th} provides necessary conditions for a stable map in $\ov\cM_{g,k}(\cZ_0,A)$ to be the Gromov(-type) limit of a sequence of stable maps in~(\ref{Family*_e}). We claim that these conditions are also virtually sufficient.
We describe the deformation-obstruction long exact sequence in Section~\ref{Deformation_s} and show that the moduli space $\ov\cM^{\tn{log}}_{g,k}(\cZ_0,A)$ is of the expected dimension equal to (\ref{exp-dim_e}).  
The claim is that $\ov\cM^{\tn{log}}_{g,k}(\cZ_0,A)$ admits a VFC that is cobordant to VFC of smooth fibers in (\ref{Family*_e}). Assuming that, we provide an explicit formula for the contributions of the virtually main components of $\ov\cM^{\tn{log}}_{g,k}(\cZ_0,A)$ to its VFC; see Formula~(\ref{Formula_e}). In Section~\ref{Cubic_s}, we work out the details for the same non-trivial example considered in \cite{ACGS} to highlight the similarities and differences. Constructing VFC and proving the degeneration formula (\ref{Formula_e}) needs a gluing theorem (with the space of gluing parameters described in (\ref{GluignEquation_e2})) that will appear in a future work. We will also need to introduce a generalized version of Kuranishi structures/space that allows ``toroidal singularities".
%------------------------------------------------------------------------------------------------------
\section{Analytical log moduli spaces}\label{LogSCVar_s}
%------------------------------------------------------------------------------------------------------
In this section, associated to any $d$-semistable SNC K\"ahler  variety $X_\eset\!=\!\bigcup_{i\in \cI} X_i$, $g,k\!\in\!\N$, and $A\!\in\!H_2(X_\eset,\Z)$, we construct the (analytic) moduli space $\ov\cM_{g,k}^{\tn{log}}(X_\eset,A)$ of $k$-marked genus $g$  degree $A$ log holomorphic curves (as a set).\\

\noindent
Given an SNC variety $X_\eset\!=\!\bigcup_{i\in \cI} X_i$, let 
$$
X_\partial = \bigcup_{\substack{ i,j\in \cI\\ i\neq j}} X_{ij}
$$
denote its singular locus.
In \cite{Fr}, associated to any SNC variety $X_\eset$, Friedman constructs a holomorphic line bundle 
$$
\cO_{X_\partial}(X_\eset) \lra X_{\partial}
$$
such that 
\bEq{PofO_e}
\cO_{X_\partial}(X_\eset)|_{X_I}= \bigotimes_{i\in I} \cN_{X_{I-i}}X_I \otimes \bigotimes_{i\in \cI-I} \cO_{X_I}(X_{I+i})\quad \forall~I\!\subset\! \cI,~|I|\!\geq\!2,
\eEq
where $\cO_{X_I}(X_{I+i})$ is the line bundle associated to the smooth divisor $X_{I+i}\!\defeq\! X_{I\cup \{i\}}\!\subset\! X_I$. \\

\noindent
If $X_\eset\!=\!\cZ_0$ is the central fiber of a smoothing $\cZ$ as in (\ref{pi_e}), $\cO_{X_\partial}(X_\eset)$ is the restriction to $X_\partial$ of $\cO_{\cZ}(\cZ_0)$ and  (\ref{LBC_e}) coincides with (\ref{PofO_e}). An SNC variety is smoothable only if $\cO_{X_\partial}(X_\eset)$ is trivial, but the converse is not true; see \cite[Sec~3]{PP} for examples.
An SNC variety $X_\eset$ is called \textbf{d-semistable} if $\cO_{X_\partial}(X_\eset)$ is isomorphic to the trivial line bundle; see \cite[Dfn.~(1.13)]{Fr}.  
Regarding the connection between the d-semistability condition and log geometry, the result is that (see \cite[Thm.~5.9]{ACGHOSS}): \textit{if $X_\eset$ is a normal crossings variety over the spectrum of an algebraically closed field, then $X_\eset$ can be equipped with a log structure over the standard log point, such that the structure morphism is log smooth if and only if $X$ is d-semistable}. We will use a reinterpretation of this statement in Section~\ref{Deformation_s}.\\

\noindent
For $I=\{i\}$, we define
\bEq{SpecialcN_e}
\cN_{X_\eset} X_{i}\equiv \bigg(\bigotimes_{j\in \cI-\{i\}} \cO_{X_{i}}(X_{ij})\bigg)^{-1}.
\eEq
With this convention, the line bundle $\cO_{X_\partial}(X_\eset)$ extends to $X_\eset$ and the trivialization (\ref{PofO_e}) compatibly extends to the case where $|I|=1$. If a smoothing $\cZ$ of $X_\eset$ as in (\ref{pi_e}) is given, then $\cN_{X_\eset} X_{i}$ coincides with $\cN_i$ in (\ref{LBC_e}) and $\cN_{i}|_{X_I}=\cN_{X_{I-i}}X_I$ for all $i\in I$ and $|I|>1$. Therefore, for simplicity, in the following we will write $\cN_i$ instead of $\cN_{X_\eset} X_{i}$.\\

\noindent
In the following construction, for each $\eset\!\neq\! I\!\subset\!\cI$ and every $i\!\in\!\cI-I$, we need to fix a holomorphic section  $\xi_{I,i}$ of $\cO_{X_I}(X_{I+i})$ vanishing (to order $1$) along $X_{I+i}$; the section $\xi_{I,i}$ is unique up to multiplication by a constant. Because of the natural isomorphism  
$$
\cO_{X_I}(X_{I+i})|_{X_J}= \cO_{X_J}(X_{J+i})\qquad \forall~\eset\neq I\subset J, ~i\in \mc{I}-J,
$$
we choose the set $\{\xi_{I,i}\}$ so that 
\bEq{ComaptibleXi_e}
\xi_{I,i}|_{X_J}= \xi_{J,i}\qquad \forall~\eset \neq I\subset J, ~i\in \mc{I}-J.
\eEq

\bRm{Symplectic}
In \cite{FMZDiv} and \cite{FMZDiv2}, with McLean and Zinger,  we introduced topological notions of normal crossings symplectic divisor and variety and established that they are equivalent, in a suitable sense, to the desired geometric notions. In \cite{FMZSum}, we showed that the direct analogue of d-semistability condition is the only obstruction to smoothability in the symplectic topology category. The process of constructing a 1-parameter family of smoothings $\cZ$ in \cite{FMZSum} is a multifold analogue of the now classical (two-fold) symplectic sum construction. Conversely, we introduced a multifold symplectic cut construction in \cite{FZCut} that, given certain configuration of Hamiltonian torus actions, degenerates a smooth target into an SNC symplectic variety. Subject to the existence of an appropriate\footnote{We need an almost complex structure $J$ on $\cZ$ such that the projection map $\pi$ in (\ref{pi_e}) is $(\mfi,J)$-holomorphic and the Nijenhueis tensor of $J$ vanishes to the first order (at least) along $\cZ_0\!\subset\cZ$; see \cite[(1.3)]{FRelative}.} class of almost complex structures on $\cZ$, the results of this paper and the rest of the claims will extend to the symplectic category; see \cite{FTowards} for an outline.
\eRm

\noindent
Let $\Gamma\!=\!\Gamma(\V,\E,\L)$ be a graph with the set of vertices $\V$, edges $\E$, and legs $\L$; the latter, also called flags or roots, are half edges that have a vertex at one end and are open at the other end. Let $\uvec{\E}$ be the set of edges with an orientation. Given an oriented edge $\uvec{e}\!\in\!\uvec{\E}$, let $\ucev{e}$ denote the same edge $e$ with the opposite orientation. For each $\uvec{e}\!\in\!\uvec{\E}$, let $v_1(\uvec{e})$ and $v_2(\uvec{e})$ in $\V$ denote the starting and ending points of the arrow, respectively. For $v,v'\!\in\!\V$, let $\E_{v,v'}$ denote the subset of edges between the two vertices and $\uvec{\E}_{v,v'}$ denote the subset of oriented edges from $v$ to $v'$. For every $v\!\in\!\V$, let $\uvec{\E}_{v}$ denote the subset of oriented edges starting from $v$.\\

\noindent
A \textbf{genus labeling} of $\Gamma$ is a function $g\colon\!\V\!\lra\! \N$. An \textbf{ordering of the legs} of $\Gamma$ is a bijection $a\colon \!\L\!\lra\! \{1,\ldots, |\L|\}$.
If a decorated graph $\Gamma$ is connected, the \textbf{arithmetic} genus of $\Gamma$ is 
$$
g=g_\Gamma= \sum_{v\in \V} g_v\! +\! \tn{rank} ~H_1(\Gamma,\Z),
$$
where $H_1(\Gamma,\Z)$ is the first homology group of the underlying topological space of $\Gamma$. Figure~\ref{labled-graph_fg}-left illustrates a labeled graph with $2$ legs. \\

\begin{figure}
\begin{pspicture}(8,1.3)(11,3)
\psset{unit=.3cm}

\pscircle*(35,6){.25}\pscircle*(39,6){.25}
\pscircle*(35,9){.25}\pscircle*(41.5,9){.25}\pscircle*(39,9){.25}
\psline[linewidth=.05](39,6)(35,9)
\psline[linewidth=.05](39,6)(39,9)\psline[linewidth=.05](39,6)(41.5,9)
\psarc[linewidth=.05](33,7.5){2.5}{-36.9}{36.9}\psarc[linewidth=.05](37,7.5){2.5}{143.1}{216.0}
\psline[linewidth=.05](41.5,9)(40.5,10.5)\psline[linewidth=.05](41.5,9)(42.5,10.5)
\rput(40.5,11.2){1}\rput(42.5,11.2){2}
\rput(35,5){$g_4$}\rput(39.5,5){$g_5$}
\rput(34.5,10){$g_1$}\rput(42.3,8.5){$g_3$}
\rput(39,10){$g_2$}

\pscircle*(55,6){.25}\pscircle*(59,6){.25}
\pscircle*(55,9){.25}\pscircle*(61.5,9){.25}\pscircle*(59,9){.25}
\psline[linewidth=.05](59,6)(55,9)
\psline[linewidth=.05](59,6)(59,9)\psline[linewidth=.05](59,6)(61.5,9)
\psarc[linewidth=.05](53,7.5){2.5}{-36.9}{36.9}\psarc[linewidth=.05](57,7.5){2.5}{143.1}{216.0}
\psline[linewidth=.05](61.5,9)(60.5,10.5)\psline[linewidth=.05](61.5,9)(62.5,10.5)
\rput(60.5,11.2){1}\rput(62.5,11.2){2}
\rput(55,5){$(g_4,A_4)$}\rput(59.5,5){$(g_5,A_5)$}
\rput(54,10){$(g_1,A_1)$}\rput(63.8,8.5){$(g_3,A_3)$}
\rput(58.3,10){$(g_2,A_2)$}
\end{pspicture}
\caption{On left, a labeled graph $\Ga$ representing elements of $\ov\cM_{g,2}$. On right, a labeled graph $\Ga$ representing elements of $\ov\cM_{g,2}(X,A)$.}
\label{labled-graph_fg}
\end{figure}

\noindent 
Such decorated graphs $\Gamma$ characterize different topological types of nodal marked curves 
$$
(\Si,\vec{z}\!=\!(z_1,\ldots,z_k))
$$ 
in the following way.
Each vertex $v\!\in\!\V$ corresponds to a smooth\footnote{We mean a smooth closed oriented surface.} component $\Si_v$ of $\Si$ with genus $g_v$. Each edge $e\!\in\!\E$ corresponds to a node $q_e$ obtained by connecting $\Sigma_v$ and $\Sigma_{v'}$ at the points $q_{\uvec{e}}\!\in\!\Si_v$ and $q_{\scz\ucev{e}}\!\in\!\Si_{v'}$, where $e\!\in\!\E_{v,v'}$ and $\uvec{e}$ is an orientation on $e$ with $v_1(\uvec{e})\!=\!v$. The last condition uniquely specifies $\uvec{e}$ unless $e$ is a loop connecting $v$ to itself. Finally, each leg $l\!\in\!\L$ connected to the vertex $v_l$ corresponds to a marked point $z_{a_l}\!\in \Sigma_{v_l}$ disjoint from the connecting nodes. If $\Si$ is connected, then $g_\Gamma$ is the arithmetic genus of $\Si$. Thus we have
\bEq{nodalcurve_e}
(\Si,\vec{z})\! =\! \coprod_{v\in \V}(\Si_v,\vec{z}_v,{q}_v)/\sim, \quad  q_{\uvec{e}}\!\sim\! q_{\scz\ucev{e}}\quad\forall~e\!\in\!\E,
\eEq
where 
$$
\vec{z}_v\!=\!\vec{z}\cap \Sigma_v\quad\tn{and}\quad \quad {q}_v=\{q_{\uvec{e}}\colon \uvec{e}\!\in\!\uvec{\E}_v\}\qquad \forall~v\!\in\!\V
$$
are the set of marked and nodal points on $\Si_v$, respectively.
In this situation, we say $\Gamma$ is the \textbf{dual graph} of $(\Si,\vec{z})$.
We treat $q_v$ as an un-ordered set of marked points on $\Si_v$. If we fix an ordering on the set $q_v$, we denote the ordered set by $\vec{q}_v$. Figure~\ref{labled-curve_fg} illustrates a nodal curve with $(g_1,g_2,g_3,g_4,g_5)=(0,2,0,1,0)$ corresponding to Figure~\ref{labled-graph_fg}-left. \\

\noindent 
Similarly, for nodal marked curves mapping into a topological space $X$, we consider similar decorated graphs where the vertices carry an additional \textbf{degree labeling} 
$$
A\colon \V\lra H_2(X,\Z),\quad v\lra A_v,
$$
recording the homology class of the image of the corresponding component. Figure~\ref{labled-graph_fg}-right illustrates a dual graph associated to a marked nodal map over the graph on the left. \\

\begin{figure}
\begin{pspicture}(8,1.7)(11,4.5)
\psset{unit=.3cm}

\pscircle*(52.4,8.2){.25}\pscircle*(52.4,11.8){.25}
\pscircle(50,10){3}\psellipse[linestyle=dashed,dash=1pt](50,10)(3,1)
\pscircle(54,7){2}\psellipse[linestyle=dashed,dash=1pt](54,7)(2,.77)
\pscircle*(55,8.6){.25}\rput(56,9){$z^1$}
\pscircle*(55,5.4){.25}\rput(56,5.4){$z^2$}

\pscircle(45,10){2}\psellipse[linestyle=dashed,dash=1pt](45,10)(2,.77)
\pscircle*(47,10){.25}\pscircle*(43.4,8.8){.25}\pscircle*(43.4,11.2){.25}
\psarc(41.8,7.6){2}{270}{71}
\psarc(41.8,12.4){2}{-71}{90}
\psarc(41.8,10){4.4}{90}{270}\psarc(42.8,10){.6}{120}{240}
\psarc(40,8.7){1.5}{60}{120}
\psarc(40,11.2){1.5}{225}{315}

\psarc(54,13){2}{60}{300}
\psarc(58,13){2}{240}{120}
\psarc(56,16.464){2}{240}{300}
\psarc(56,9.546){2}{60}{120}
\psarc(54,11.7){1.5}{60}{120}
\psarc(54,14.2){1.5}{225}{315}
\psarc(58,11.7){1.5}{60}{120}
\psarc(58,14.2){1.5}{225}{315}

\end{pspicture}
\caption{A nodal curve in $\ov\cM_{4,2}$.}
\label{labled-curve_fg}
\end{figure}

\noindent
Assume $X_{\eset}\!=\!\bigcup_{i\in \cI} X_i$ is an SNC variety and $\Si$ is an irreducible smooth curve. Then every holomorphic map $u\colon\! \Si\!\lra\! X_\eset$ has a well-defined \textbf{depth} $\eset\!\neq\! I\!\subset\!\cI$, which is the maximal subset of $\cI$ such that $\tn{Image}(u)\!\subset\!X_I$. Similarly, we say a point $x\!\in\!\Si$ has \textbf{depth $I$}, if $X_I$ is the minimal stratum containing $u(x)$. Let $\cP^*(\cI)$ be the set of non-empty subsets of $\cI$. The dual graph of $(u,\Si)$ carries additional labelings
$$
I\colon \V,\E\lra \cP^*(\cI), \qquad v\lra I_v\quad \forall v\!\in\!\V, \qquad e\lra I_e\quad \forall e\!\in\!\E,
$$
recording the depths of smooth components and nodes of $\Si$. \\

\noindent
Given a finite set $\cI$ and a ring $R$, let
$$
R^\cI_\bullet=\big\{r=(r_j)_{j\in \cI}\colon \sum_{j\in \cI} r_j=0\big\}\subset R^\cI.
$$
For every $i\!\in\!\cI$, the projection
\bEq{pibullet_e}
\pi_i\colon R^\cI_\bullet\lra R^{\cI-\{i\}},\quad (r_j)_{j\in \cI} \lra (r_j)_{j\in \cI-\{i\}}
\eEq
is an isomorphism. For every subset $\cI'\!\subset\!\cI$, the natural injective homomorphism $R^{\cI'}\!\hookrightarrow\!R^{\cI}$ restricts to a homomorphism
$R_\bullet^{\cI'}\!\hookrightarrow\!R^{\cI}_\bullet$. Therefore, via this inclusion, $R_\bullet^{\cI'}$ can be thought of as a subspace of $R^{\cI}_\bullet$.  We will use the projection maps in (\ref{pibullet_e}) to identify each component of a pre-log moduli space in $X_\eset$ with a fiber  product of relative spaces in $\{X_i\}_{i\in \mc{I}}$. \\

\noindent
Assume $X_\eset\!=\!\bigcup_{i\in \cI} X_i$ is a d-semistable SNC K\"ahler variety and fix a trivialization of $\cO_{X_\partial}(X_\eset)$.
Let $u\colon\!\Si\!\lra\!X_\eset$ be a holomorphic map of depth $I\!\subset\!\cI$ with smooth domain. Then, for every $i\!\in\!\cI\!-\!I$, the function
\bEq{Ordx_e1}
\ord_{u}^i\colon\Si\lra \N, \quad \ord_{u}^i(x)=\tn{ord}_x(u,X_i),
\eEq
recording the contact order of $u$ with the smooth divisor $X_{I+i}\!\subset\! X_I$ at $x$ is well-defined. It coincides with the order $\tn{ord}_x(u^* \xi_{I;i})$ of $\xi_{I;i}$ at $x$ where $\xi_{I;i}$ are the sections in (\ref{ComaptibleXi_e}).  \\

\noindent
For every $\eset \!\neq\!I\!\subset\!\cI$, each $i\!\in\!I$, and a given meromorphic section $\ze$ of  $u^*\cN_{X_{I-i}} X_I$, we will also need the function 
\bEq{Ordx_e2}
\ord_{\ze}\colon\Si\lra \Z, \quad \ord_{\ze}(x)= \tn{ord}_x(\ze),
\eEq
recording the vanishing order of $\ze$ at $x$ (which is negative if $\ze$ has a pole at $x$).\\

\noindent
Given a holomorphic line bundle $\cL\!\lra\!\Si$, let $\Om_{\tn{mero}}(\Si,\cL)$ denote the space of non-trivial meromorphic sections of $\cL$.

\noindent
\bDf{LogTuple_dfn}
For each $\eset\!\neq\!I\!\subset\!\cI$, a \textbf{log holomorphic tuple} $(u,\ze, \Si, q)$ of depth $I$ consists of a smooth irreducible curve $\Si$,  a finite set of distinct points $q\!=\!\{q_1,\ldots,q_\ell\}$ on $\Si$,  a depth $I$ holomorphic map $u\colon\!\Si\!\lra\!X_I$, and a set of meromorphic sections
$$
\ze \equiv \big(\ze_i\big)_{i\in I}\in    \prod_{i\in I} \Om_{\tn{mero}}(\Si,u^*\cN_{X_{I-i}}X_I)
$$
such that:\newline 
(a) via the identification (\ref{PofO_e}) and the given trivialization of $\cO_{X_\partial}(X_\eset)$,  we have
$$
\bigotimes_{i\in I} \ze_i \otimes \bigotimes _{j \in \cI-I} u^*\xi_{I,j}=u^*(\la_{\cO_{X_I}})
$$ 
for some $\la\!\in\! \C^*$,\\
(b) for all $x\!\in\!\Si$
\bEq{Intersection-Condition_e}
\tn{if}~~\ord_{u,\ze}(x)\!\neq \!0~~\Rightarrow~~x\!\in\!q,
\eEq
where the vector-valued \textbf{order function}
$$
\ord_{{u,\ze}}(x)\!=\!\big((\ord_{u}^j(x))_{j\in \cI-I},(\ord_{\ze_i}(x))_{i\in I}\big)=\!\big((\ord_{u^*\xi_{I;j}}(x))_{j\in \cI-I},(\ord_{\ze_i}(x))_{i\in I}\big) \!\in\!\Z^\cI_\bullet 
$$ is defined via (\ref{Ordx_e1}) and (\ref{Ordx_e2}).
\eDf

\noindent
\bRm{elboration_rmk}
Let us elaborate on Definition~\ref{LogTuple_dfn}.
\bIt
\item In (b), the reason we get a vector in $\Z^\cI_\bullet$  is Condition (a).
\item If $X_\eset=\cZ_0$ is the central fiber of a semistable degeneration and $u$ is of degree $A\!\in\!H_2(\cZ_0,\Z)$, then (\ref{Intersection-Condition_e}) implies that
$$
(A\cdot X_i)_{i\in \cI}=\sum_{q_a\in q} \ord_{u,\ze}(q_a) \in \Z^{\cI}_\bullet.
$$
Note that $A\cdot \cZ_0\!=\! \sum_{i\in \cI} A\cdot X_i \!=\!0$.
\item Changing any of the meromorphic sections $\ze_i$ with a non-zero constant multiple of that has no effect on any of the conditions. 
\item If $I\!=\!\{i\}$, by Condition~(a) and (\ref{SpecialcN_e}), the section $\ze_i$ is uniquely determined by $u$ and $\la$. For $I$ with $|I|\!>\!1$, specifying $|I|\!-\!1$ of sections $(\ze_i)_{i\in I}$ and $\la$ will uniquely determine the remaining one. So there is a redundancy in Definition~\ref{LogTuple_dfn}. The reason for this redundancy is to avoid a non-symmetric definition that depends on the choice of $i\!\in\!I$.
\item For a fixed set $q=\{q_1,\ldots,q_\ell\}\!\subset \!\Si$, fixed $u$, and a fixed set of order vectors $s_1,\ldots,s_\ell\in \Z^\cI_\bullet$, up to multiplication by non-zero constants, there is at most one set of meromorphic sections $\ze$ such that 
$$
\ord_{u,\ze}(q_a)=s_a\quad \forall~a=1,\ldots,\ell\qquad \tn{and}\qquad \ord_{u,\ze}(x)=0\quad \forall x\in \Si-q;
$$
see \cite[Lmm.~2.15]{FRelative}.
Therefore, the analytical logarithmic moduli spaces below can be described without mentioning the meromorphic sections $\ze$. Since these sections appear in several calculations, it is convenient to consider $\ze$ as part of the defining data of a log map.
\item The set of vectors
$$
\mfs\!=\!\big\{s_i=(s_{ij})_{j\in \cI}\big\}_{1\leq i\leq \ell}\!\subset\!(\Z^\cI_\bullet)^\ell 
$$
describe the ``contact type" of the log holomorphic tuple $(u,\ze, \Si, q)$ with the singular locus $X_\partial$ and will play an important role in defining log moduli spaces. They also appear in the relative compactification of \cite{FRelative} at the marked points. Here, they only appear at the nodal points.
\eIt
\eRm

\bDf{PreLogMap_dfn} 
Let $X_\eset\!=\!\bigcup_{i\in \cI} X_i$ be a d-semistable SNC K\"ahler variety and fix a trivialization of $\cO_{X_\partial}(X_\eset)$. Suppose 
$$
C\!\equiv\!(\Si,\vec{z})=\bigg(\coprod_{v\in \V} C_v\equiv(\Si_v,\vec{z}_v,q_v)\bigg)/\sim,\quad   q_{\uvec{e}}\!\sim\! q_{\scz\ucev{e}}\quad\forall~\uvec{e}\!\in\!\uvec{\E},
$$ 
is a $k$-marked connected nodal  curve with smooth components $C_v$ and dual graph $\Gamma\!=\!\Gamma(\V,\E,\L)$ as in (\ref{nodalcurve_e}). A \textbf{pre-log map} from $C$ to $X_\eset$ is a collection 
\bEq{fplogSetUp_e}
f\equiv \big(f_{v}\!\equiv\! (u_v,\ze_v,C_v)\big)_{v\in \V}
\eEq
such that  
\bEn
\item for each $v\!\in\!\V$, $(u_v,\ze_v=(\ze_{v,i})_{i\in I_v},\Si_v, q_v)$ is a depth $I_v$ log holomorphic tuple as in Definition~\ref{LogTuple_dfn} for some fixed\footnote{i.e., independent of the choice of $v\in \V$.} $\la=\la(f)$,
\item\label{MatchingPoints_l} $u_v(q_{\uvec{e}})\!=\!u_{v'}(q_{\scz\ucev{e}})$ for all $\uvec{e}\!\in\! \uvec{\E}_{v,v'}$;
\item\label{MatchingOrders_l} $s_{\uvec{e}}\!\equiv\! \ord_{u_v,\ze_v}(q_{\uvec{e}})\!=\! -\ord_{u_{v'},\ze_{v'}}(q_{\scz\ucev{e}})\!\equiv\!-s_{\scz\ucev{e}}$ for all $v,v'\!\in\!\V$ and $\uvec{e}\!\in\! \uvec{\E}_{v,v'}$;
\eEn
\eDf
\noindent
In other words, a pre-log map is a nodal holomorphic map into $X_{\eset}$ with a bunch of meromorphic sections (satisfying Definition~\ref{LogTuple_dfn}.(a) on each smooth component), opposite contact orders at the nodes, and trivial contact orders at every other point (including the marked points).

\bRm{IevsIv_rmk}
For every $v\!\in\!\V$ and $\uvec{e}\!\in\!\uvec{\E}_v$, let
\bEq{seDf_e}
s_{\uvec{e}}= (s_{\uvec{e},i})_{i\in \cI}=\big((\ord_{u_v}^i(q_{\uvec{e}}))_{i\in \cI-I_v},(\ord_{\ze_{v,i}}(q_{\uvec{e}}))_{i\in I_v}\big) \!\in\!\Z^{\cI}_\bullet
\eEq
be the contact order vector at the nodal point $q_{\uvec{e}}\!\in\!\Si_v$.
For $\uvec{e}\!\in\! \uvec{\E}_{v,v'}$, since $u_v$ and $u_{v'}$ have image in $X_{I_{v}}$ and $X_{I_{v'}}$, respectively, by Condition~\ref{MatchingPoints_l} above, we have
$$
u(q_e)\!=\!u_v(q_{\uvec{e}})\!=\!u_{v'}(q_{\scz\ucev{e}})\in X_{I_{v}}\cap X_{I_{v'}}= X_{I_{v}\cup I_{v'}};
$$
i.e. $I_{e}\!\supset\!I_{v}\!\cup\! I_{v'}$. If $i\!\in\!\cI\!\setminus\!I_{v}\!\cup\! I_{v'}$, by (\ref{Ordx_e1}), we have 
$$
s_{\uvec{e},i},s_{\scz\ucev{e},i}\geq 0.
$$
Therefore, by Condition~\ref{MatchingOrders_l} above, they are both zero, i.e.
 $$ 
 I_{e}\!=\!I_{v}\!\cup\! I_{v'}\quad\tn{and}\quad s_{\uvec{e}}\in \Z^{I_{e}}_\bullet  \subset \Z^{\cI}_\bullet\qquad \forall~\uvec{e}\!\in\! \uvec{\E}_{v,v'}.
$$
 \eRm
 
 \noindent
Two pre-log maps $(u,\ze,C)\!\equiv\! \big(u_v,\ze_v,C_v\big)_{v\in \V}$ and $(\wt{u},\wt{\ze},\wt{C})\!\equiv\!\big(\wt{u}_v,\wt{\ze}_v,\wt{C}_v\big)_{v\in \V}$  with isomorphic decorated dual graphs $\Gamma$ are equivalent if there exists a biholomorphic identification 
\bEq{NodalEquivalence_e}
\big(h\colon \wt{C}\lra C\big)\equiv \big(h_v\colon \wt\Si_v\!\lra\!\Si_{h(v)}\big)_{v\in \V}
\eEq \vskip-.1in
\noindent
such that
$$
h(\wt{z}_a)\!=\!z_a\quad \forall a\!=\!1,\ldots,k,\quad u\circ h\!=\!\wt{u},\quad h_v^*\ze_{h(v),i}= c_{v,i} \wt{\ze}_{v,i}\quad \forall~v\!\in\!\V,~i\!\in\!I_v,~\tn{some}~c_{v,i}\in \C^*.
$$
A pre-log map $f$ is \textbf{stable} if the group of self-equivalences $\aut(f)$ is finite. 
By the fifth bullet in Remark~\ref{elboration_rmk}, a pre-log map is stable if and only if the underlying marked nodal map is stable. Clearly, the automorphism group of a pre-log map is a subgroup of the automorphism group of the underlying marked nodal map.
The equivalence class of a pre-log map is called a \textbf{pre-log curve}. 
For every choice of a decorated dual graph $\Gamma$, we denote the space of $k$-marked degree $A$ pre-log holomorphic curves of type $\Gamma$ by 
$$
\cM^{\tn{plog}}_{g,k}(X_\eset,A)_\Gamma.
$$

\noindent
A basic dimension counting reveals that the expected dimension of $\cM^{\tn{plog}}_{g,k}(X_\eset,A)_\Gamma$ could be much larger than (\ref{exp-dim_e}). In Definition~\ref{LogMapVar_dfn} below, we take out a subspace that would give us a nice compactification with the correct expected dimension. \\

\noindent
Let  
$$
\Big(\bigoplus_{v\in \V}\Z^{I_v}\Big)^*=\Big\{(s_v)_{v\in \V}\in \bigoplus_{v\in \V}\Z^{I_v}\colon \sum_{i\in I_v} s_{v,i}=\sum_{i\in I_{v'}} s_{v',i}\quad \forall v,v'\in \V\Big\},
$$
and 
\bEq{ToricPi_e}
\pi\colon \Big(\bigoplus_{v\in \V}\Z^{I_v}\Big)^*\lra \Z, \qquad (s_v)_{v\in \V}\lra \sum_{i\in I_v} s_{v,i}\,.
\eEq
Associated to the decorated dual graph $\Gamma\!=\!\Gamma(\V,\E,\L)$ of a pre-log map and an arbitrary orientation $O\!\equiv\!\{\uvec{e}\}_{e\in \E} \subset \uvec{\E}$ on the edges, we define a $\Z$-linear map 
\bEq{vrStar_e}
\vr \colon \D\defeq \Z^\E\oplus \Big(\bigoplus_{v\in \V}\Z^{I_v}\Big)^*\lra \T\defeq \bigoplus_{e\in \E} \Z_\bullet^{I_e}
\eEq
in the following way: 
\bIt
\item $\vr$ maps the generator $1_{e}$ of the $e$-th component in the summand $\Z^\E\!\subset\! \D$ to the contact vector $s_{\uvec{e}}\in \Z_\bullet^{I_e}\!\subset\! \T$, where $\uvec{e}$ is the chosen orientation on $e$ in $O$;
\item  if $\uvec{e}$ is the chosen orientation on $e$ in $O$, the $e$-th component of $\vr(\oplus_{v\in \V} s_v)$ is  
$$
s_{v_1(\uvec{e})}-s_{v_2(\uvec{e})}\in \Z_\bullet^{I_e}.
$$
\eIt

\noindent
Let 
$$
\Lambda\!=\!\tn{image}(\vr),\quad \tn{K}\!=\tn{Ker}(\vr)\quad\tn{and}\quad  
\tn{CK}\!=\!\T/\Lambda=\tn{Coker}(\vr).
$$
The $\Z$-modules $\Lambda$, $\tn{K}$, and $\tn{CK}$ are independent of the choice of the orientation $O$ on $\E$ and are invariants of the decorated graph $\Gamma$. In particular,
\bEq{kernel_e3}
\tn{K}=\lrc{\big((\la_e)_{e\in \E},(s_v)_{v\in \V}\big)\!\in\!\ \Z^\E\oplus \Big(\bigoplus_{v\in \V}\Z^{I_v}\Big)^* \colon s_{v'}\!-\!s_{v}\!=\!\la_es_{\uvec{e}}\quad\forall~v,v'\!\in\!\V,~\uvec{e}\!\in\!\uvec{\E}_{v,v'}}.
\eEq
Replacing $\Z$ with another ring $R$ in the equations above, we denote the corresponding terms in (\ref{vrStar_e}) by $\vr_R\colon \D_R\lra \T_R$ and image/kernel/cokernel by $\Lambda_{R}$, $\tn{K}_{R}$, and $\tn{CK}_{R}$, respectively. Via the exponentiation map, let 
$$
\tn{exp}(\Lambda_{\C}) \subset  \prod_{e\in \E} (\C^*)^{I_e}_\bullet,\quad\tn{with}~~ (\C^*)^{I}_\bullet=\{(t_i)_{i\in I}\!\in (\C^*)^I: ~~\prod_{i\in I} t_i =1\},
$$
be the subgroup corresponding to the sub-Lie algebra $\La_{\C}\!\subset\!\T_{\C}$, and denote the quotient group by 
$$
\mc{G}=\!\mc{G}(\Gamma)\!=\!\tn{exp}(\tn{CK}_{\C}).
$$
In other words,
$$
\mc{G}=\frac{\prod_{e\in \E} (\C^*)_\bullet^{I_e}}{\tn{exp}(\vr)\Big((\C^*)^\E\times \big( \prod_{v\in \V} (\C^*)^{I_v}\big)^* \Big)}, 
$$
where 
$$
 \Big( \prod_{v\in \V} (\C^*)^{I_v}\Big)^*\defeq \Big\{(t_v)_{v\in \V}\in \prod_{v\in \V}(\C^*)^{I_v}\colon \prod_{i\in I_v} t_{v,i}=\prod_{i\in I_{v'}} t_{v',i}\quad \forall v,v'\in \V\Big\}.
 $$

\noindent
Below, similarly to \cite{FRelative}, to every pre-log $f$ as in Definition~\ref{PreLogMap_dfn} we associate a group element
$$
\tn{ob}_{\Gamma}(f) \in \mc{G}(\Gamma)
$$
that only depends on the equivalence class of $f$.\\

\noindent
Given a pre-log map $f\equiv \big(f_{v}\!\equiv\! (u_v,\ze_v,C_v)\big)_{v\in \V}$ as in (\ref{fplogSetUp_e}),
for each $v\!\in\!\V$ and $\uvec{e}\!\in\!\uvec{\E}_v$, let $z_{\uvec{e}}$ be an arbitrary holomorphic coordinate in a sufficiently small disk $\De_{\uvec{e}}$ around the nodal point $(z_{\uvec{e}}\!=\!0)\!=\!q_{\uvec{e}}\!\in\!\Si_v$.  

\bEn
\item By (\ref{Ordx_e2}), for every  $v\!\in\!\V$, $\uvec{e}\!\in\!\uvec{\E}_v$, and $i\!\in\!I_v$, in a local holomorphic trivialization 
$$
u_v^*\cN_{X_{I_v-\{i\}}}X_{I_v}|_{\De_{\uvec{e}}}\approx \cN_{X_{I_v-i}}X_{I_v}|_{u(q_e)}\!\times\!\De_{\uvec{e}},
$$ 
we have 
\bEq{LocalCoord_e}
\ze_{v,i}(z_{\uvec{e}})= z_{\uvec{e}}^{s_{\uvec{e},i}} \wt{\ze}_{v,i}(z_{\uvec{e}})
\eEq
such that 
$$
0\!\neq\!  \wt{\ze}_{v,i}(0) \! \equiv\!\eta_{\uvec{e},i}\!\in\!\cN_{X_{I_v-i}}X_{I_v}|_{u(q_e)}=\cN_i|_{u(q_e)}
$$ 
is independent of the choice of  the trivialization. 
\item Similarly, 
for every  $v\!\in\!\V$, $\uvec{e}\!\in\!\uvec{\E}_v$, and $i\!\in\!I_e\!-\!I_v$,
with $\xi_{I_v,i}$ as in Definition~\ref{LogTuple_dfn}.(a),  we have
\bEq{LocalCoord_e2}
u_v^*\xi_{I_v,i}(z_{\uvec{e}})= z_{\uvec{e}}^{s_{\uvec{e},i}} \wt{\xi}_{v,i}(z_{\uvec{e}})
\eEq
such that
$$
0\!\neq\!  \wt{\xi}_{v,i}(0) \! \equiv\!\eta_{\uvec{e},i}\!\in\!\cN_{X_{I_v}}X_{I_v+\{i\}}|_{u(q_e)}=\cN_i|_{u(q_e)}.
$$  
 Note that the map $u_v$ 
has a well-defined $s_{\uvec{e},i}$-th derivative 
\bEq{normalDerivative}
\eta'_{\uvec{e},i}\!\in\!\cN_{X_{I_v}}X_{I_v+i}|_{u(q_e)}
\eEq
(with respect to the coordinate $z_{\uvec{e}}$) in the normal direction to $X_{I_v+i}\subset X_{I_v}$ at the nodal marked point $q_{\uvec{e}}$. The vector $\eta'_{\uvec{e},i}$ is a multiple of $\eta_{\uvec{e},i}$ by a factor that only depends on the choice of $\xi_{I_v,i}$.
\item Finally, for every  $v\!\in\!\V$, $\uvec{e}\!\in\!\uvec{\E}_v$, and $i\!\in\!\cI\!-\!I_e$, let $\eta_{\uvec{e},i}\!\in\!\cN_{X_{I_v}}X_{I_v+\{i\}}|_{u(q_e)}=\cN_i|_{u(q_e)}$ be the non-zero
 vector $(u_{v}^*\xi_{I_v,i})(q_{\uvec{e}})$. 
\eEn
\noindent
For each $e\!\in\!\E$, if $\uvec{e}$ is the choice of orientation on $e$ in $O$, let 
\bEq{etae_e}
\eta_{e}=\big(\eta_{\uvec{e},i}/\eta_{\scz\ucev{e},i}\big)_{i\in I_e}\!\in\!(\C^*)^{I_e}.
\eEq
By Definition~\ref{PreLogMap_dfn}.(1)  and because $s_{\uvec{e}}\in \Z^{I_e}_\bullet\subset \Z^{\mc{I}}_\bullet$ we have 
$$
\otimes_{i\in \mc{I}}~ \eta_{\uvec{e},i} = \la_{\cO_{X_\partial}(X_\eset)}|_{u(q_e)}.
$$
By (\ref{ComaptibleXi_e}), we have 
$$
\eta_{\uvec{e},i}=\eta_{\scz\ucev{e},i}\qquad \forall i\!\in\!\mc{I}-I_e
$$
Therefore, 
$$
\eta_{e}\in (\C^*)^{I_e}_\bullet \qquad \forall~e\!\in\!\E.
$$

\noindent
The tuples $\eta_e$ give rise to an element 
$$
\eta\!\equiv \!(\eta_e)_{e\in \E}\!\in\! \prod_{e\in \E} (\C^*)^{I_e}_{\bullet}.
$$ 
The action of the subgroup $\tn{exp}(\vr)\Big((\C^*)^\E\times \big(\prod_{v\in \V} (\C^*)^{I_v}\big)^* \Big)$ on $\eta$ corresponds to admissible rescalings of the sections $\ze_{v,i}$ and change of coordinates in $z_{\uvec{e}}$; i.e the class 
\bEq{PLtoG_e}
\ob_{\Gamma}(f)\!=\![\eta]
\eEq
of $\eta$ in $\mc{G}$ is independent of the choices involved. If $f$ and $f'$ are equivalent with respect to a reparametrization $h\colon\!\Si'\!\lra\!\Si$ as in (\ref{NodalEquivalence_e}), the associated group elements $\eta$ and $\eta'$, respectively, would be the same with respect to any $h$-symmetric choice of holomorphic coordinates $\{z_{\uvec{e}}\}_{\uvec{e}\in\uvec{\E}}$.
Therefore, (\ref{PLtoG_e}) only depends on the equivalence class $[f]$ of $f$ and thus yields a well-defined function 
\bEq{ObBullet_e}
\tn{ob}_{\Gamma}\colon \cM_{g,k}^{\tn{plog}}(X_\eset,A)\lra \mc{G}(\Gamma).
\eEq
\bRm{RI_rmk}
It is clear from the definition that $\ob_{\Gamma}([f])=1$ if and only if there exists a representative 
$f\equiv \big(f_{v}\!\equiv\! (u_v,\ze_v,C_v)\big)_{v\in \V}$ as in (\ref{fplogSetUp_e}),
and a choice of local coordinates $z_{\uvec{e}}$ around the nodal point $(z_{\uvec{e}}\!=\!0)\!=\!q_{\uvec{e}}\!\in\!\Si_v$ such that 
$$
\eta_{\uvec{e},i}=\eta_{\scz\ucev{e},i}\qquad \forall~e\!\in\!\E,~i\!\in\!I_e.
$$
\eRm

\bDf{LogMapVar_dfn}
Let $X_\eset\!=\!\bigcup_{i\in \cI} X_i$ be a $d$-semistable SNC variety.
 A \textbf{log} map is a pre-log  map  $f$ with the decorated dual graph $\Gamma$ such that 
\bIt
\item\label{Tropical_l2} (C1) there exist functions 
$$
s \colon\! \V\!\lra\!\R^\cI, \quad v\!\lra\!s_v,\qquad\tn{and}\qquad \la \colon\! \E\!\lra\!\R_+, \quad e\!\lra\!\la_e,
$$
such that 
\bEnalph
\item $s_v\!\in\!\R_{+}^{I_v}\!\times\!\{0\}^{\cI-I_v}$ for all $v\!\in\!\V$,
\item\label{Direction_l2} $s_{v_{2}(\uvec{e})}\!-\!s_{v_1(\uvec{e})}\!=\!\la_{e} s_{\uvec{e}}$  for every $\uvec{e}\!\in\!\uvec{\E}$;
\eEnalph  
\item\label{GObs_e2} (C2) and $\tn{ob}_{\Gamma}(f)\!=\!1\!\in\!\mc{G}(\Gamma)$.
\eIt
\eDf
\vskip.1in
\noindent
Since $\Gamma$ is connected, Condition~(C1)\ref{Direction_l2} and $s_{\uvec{e}}\!\in\!\Z^\cI_\bullet$, for all $\uvec{e}\!\in\!\uvec{\E}$, imply that 
\bEq{Constantc_e}
\sum_{i\in \cI} s_{v,i}\!=\!r
\eEq
is a  fixed positive constant $r\!\in\!\R_+$ independent of $v\!\in\!\V$. The combinatorial condition (C1) is equivalent to the basicness  condition in \cite[Dfn.~1.20]{GS}. It is interesting and important to figure out the relation between (C2) and the lift-ability conditions in \cite[Thm.~5.3.3]{ACGS}. The set of vectors $\big((s_v)_{v\in\V},(\la_{e})_{e\in \E}\big)$
satisfying (C1) is the intersection of the kernel $\tn{K}_\R$ of $\vr_\R$ with the positive quadrant in $\D_\R$.
Since this intersection is non-empty by assumption, its closure is a convex maximal rational polyhedral cone $\si(\Gamma)$ in $\tn{K}_\R$. The space of gluing parameters is a multiple of the affine toric variety $Y_\si$ with the toric fan $\si(\Gamma)\!\subset\!\tn{K}_\R$; see Section~\ref{Formula_s}. The projection map $\pi$ in (\ref{ToricPi_e}) restricts to $\si(\Gamma)\to \R_{\geq 0}$; it gives a projection map between toric varieties $Y_\si \lra \C$ ($\C$ is the toric variety associated to the toric fan $\R_{\geq 0}\subset \R$) that, via gluing, corresponds to the fiberation (\ref{pi_e}).
 \\

\noindent
We will denote the subset of log curve in $\cM^{\tn{plog}}_{g,k}(X_\eset,A)_\Gamma$ by $\cM^{\tn{log}}_{g,k}(X_\eset,A)_\Gamma$. In other words, we require $\Gamma$ to satisfy the combinatorial condition (C1) and 
$$
\cM^{\tn{log}}_{g,k}(X_\eset,A)_\Gamma=\tn{ob}_{\Gamma}^{-1}(1)\subset \cM^{\tn{plog}}_{g,k}(X_\eset,A)_\Gamma.
$$
In Section~\ref{Deformation_s}, we show that the expected $\C$-dimension of $\cM^{\log}_{g,k}(X_\eset,A)_\Gamma$ is equal to 
$$
c_1^{T^{\log}X_\eset}(A)+(\dim_\C X_\eset-3)(1-g)+k - (\dim_\R \tn{K}_{\R}-1),
$$
where $T^{\log}X_\eset$ is the log tangent bundle of $X_\eset$. Therefore, $\cM^{\log}_{g,k}(X_\eset,A)_\Gamma$ is virtually a main stratum if and only if 
$\tn{K}\!\cong\! \Z$. \\

\noindent
If $\cZ$ is a semistable smoothing of $\cZ_0\!=\!X_\eset$ as in (\ref{pi_e}), $c_1^{T^{\log}X_\eset}(A)$ coincides with $c_1^{T\cZ_\la}(A)$, for all $\la\!\in\!\De^*$. Therefore, the expected dimension of the \textbf{analytical log moduli space}
$$
\ov\cM^{\log}_{g,k}(X_\eset,A)=\bigcup_{\Gamma} \cM^{\log}_{g,k}(X_\eset,A)_\Gamma
$$ 
coincides with the expected dimension of $\ov\cM_{g,k}(\cZ_\la,A)$. Different components in this union, however, will contribute with different weights to the VFC of $\ov\cM_{g,k}(\cZ_\la,A)$. The degeneration formula (\ref{Formula_e}) describes these weights. 

\bRm{vrs_rmk}
There is a slightly different map associated to $\Gamma$ that will be useful in Section~\ref{Formula_s}. Let
\bEq{DtoT_e}
\vr_\bullet \colon \D_\bullet=\Z^\E\oplus \bigoplus_{v\in \V}\Z^{I_v}_\bullet \lra \T\equiv\bigoplus_{e\in \E} \Z_\bullet^{I_e}.
\eEq
denote the restriction of $\vr$ to $\D_\bullet\subset \D$. The following diagram commutes 
\bEq{CD_e}
\xymatrix{
\D_\bullet\ar^{\vr_\bullet}[rr]\ar[d]&& \T\ar^{\cong}[d]\\
\D\ar^{\vr}[rr]\ar[d]^{0\oplus \pi}&& \T\\
\Z\,. && \\
}
\eEq
Fixing a vector in the interior of $\si(\Gamma)$ gives us a decomposition
$$
\tn{K}_\R(\vr)\cong \tn{K}_\R(\vr_\bullet)\oplus \R
$$
such that $\si(\Gamma)$ is a cone containing the ray $0\oplus \R_{\geq 0}$.  Therefore, $\cM^{\log}_{g,k}(X_\eset,A)_\Gamma$ is virtually a main stratum if and only if 
$\tn{K}_\bullet\!=\! 0$. 
\end{remark}

%------------------------------------------------------------------------------------------------------
\section{Proof of the main theorem}\label{Proof_s}
%------------------------------------------------------------------------------------------------------

In this section, by following and adjusting the steps of the proof \cite[Thm.~1.3]{FRelative}, we prove Theorem~\ref{Compactness_th}. \\

\noindent
Throughout this section, for the cases that invlove the smoothing of $X_\eset$, let $\pi\colon \cZ \lra \De$ be a K\"ahler semistable degeneration as in (\ref{pi_e}) with compact fibers and the SNC central fiber 
$$
\cZ_0\!\defeq\!\pi^{-1}(0)=X_\eset\defeq\bigcup_{i\in \cI} X_i.
$$ 

\noindent
Assume $X$ is an $(n+1)$-dimensional smooth complex variety and $D\subset X$ is a smooth complex divisor. For any $p\in D$, let $U$ be an open set around $p$ with local  coordinates  $(x_0,\ldots,x_{n})$ given by a chart map
$$
\varphi\colon U\lra \C^{n+1}
$$ 
such that $D\cap U \equiv \varphi^{-1}(x_{0}=0)$. We call such a pair $(U,\varphi)$ a \textbf{$D$-compatible chart} around $p$. A $D$-compatible chart $(U,\varphi)$ induces a holomorphic trivialization 
\bEq{dcN-XD_e}
\nd_{\cN_XD}\varphi \colon \cN_XD|_{U\cap D}\lra \varphi(U \cap D)\times \C \subset \C^{n+1}.
\eEq

\noindent
Before we state the next proposition, we need to elaborate on the topological aspects of the Gromov convergence \cite{Gr} and setup the notation. For more details, we refer to \cite[Sec. 3.1]{FRelative}.\\

\noindent
Suppose 
\bEq{OSeq_e}
\Big(f_a\in \ov\cM_{g,k}(\cZ_{\la_a},A)\Big)_{a\in \N}~,\qquad f_a= \big(u_{a,v'},C_{a,v'}=(\Si_{v'},\vec{z}_{v'})\big)_{v'\in \V'},
\eEq
is a sequence of stable maps with a fixed decorated dual graph $\Gamma'\!=\!\Gamma(\V',\E',\L')$ that Gromov converges to the stable map
\bEq{LSeq_e}
f=\big(u_{v},C_{v}\big)_{v\in \V}\in\ov\cM_{g,k}(\cZ_0,A)
\eEq
with the decorated dual graph $\Gamma\!=\!\Gamma(\V,\E,\L)$. Then, (for $a$ sufficiently large) 
all the marked domains
$$
(\Si_a=\cup_{v'\in \V'}\Si_{a,v'},\vec{z}_a=\cup_{v'\in\V'} \vec{z}_{v'})
$$ 
are smoothly isomorphic to a fixed marked domain $(\Si',\vec{z}\,')$ and the domain $\Si=\bigcup_{v\in \V} \Si_{v}$ of $f$ is obtained by collapsing a set of embedded curves away from the marked and nodal points in $\Si'$. In other words, there is a continuous degeneration map 
$$
\gamma\colon \Si'\lra \Si
$$
that sends the marked points and nodal points of $\Si'$ to the (corresponding) marked and nodal points of $\Si$, and collapses some embedded circles $\{\gamma_e\}_{e \in \E^{\tn{cut}}}$ in $\Si'$ to the new nodal points $\{q_e\}_{e \in \E^{\tn{cut}}}$ in $\Si$. The map $\gamma$ gives us a similarly denoted map 
\bEq{gammaGraph_map}
\gamma^*\colon \Gamma \lra \Gamma'.
\eEq
We have 
$$
\E\cong \E'\cup \E^{\tn{cut}}, \qquad \L=\L', 
$$
such that $\gamma^*|_{\E'\subset \E}$ and $\gamma^*|_{\L}$ are isomorphisms and 
$$
\gamma^*\colon \E^{\tn{cut}} \cup \V\lra \V'.
$$
For every $v\in \V$, there exists a unique $v'\!\in\!\V'$ and a connected component $U_{v}$ of 
$\Si_{v'}-\{\gamma_e\}_{e\in \E^{\tn{cut}}}$ such that $\Si_{v}$ is obtained by collapsing the boundary circles of the closure of $U_{v}$.\\

\noindent
The goal is to, after possibly passing to a subsequence, find a set of meromorphic sections $\{\ze_{v,i}\}_{v\in \V}$ that enhances (\ref{LSeq_e}) to a log map $f_{\log}\in \cM_{g,k}^{\log}(\cZ_0,A)_{\Gamma}$.

\bDf{STNS1_dfn} 
With notation as in (\ref{OSeq_e}) and  (\ref{LSeq_e}), let 
 $$
 \ze_{v,i}\in \Om_{\tn{mero}}(\Si_{v},u_{v}^* \cN_{i}) \qquad \forall~v\!\in\! \V,~i\!\in\! I_{v}, $$  
be a set of meromorphic sections. We say (\ref{OSeq_e})  log-Gromov converges to 
 \bEq{LogLimit_e11}
f_{\log}\defeq\big(u_{v},\ze_{v}=(\ze_{v,i})_{i\in I_{v}},C_{v}\big)_{v\in \V}
\eEq
if for each $v\!\in\!\V$ and $i\!\in\! I_{v}$, with $v'\!=\!\gamma^*(v)\!\in\!\V'$, the sequence $(u_{a,v'})_{a \in \N}$ is asymptotic to $\ze_{v,i}$ in the normal direction to $X_i$ in the following sense. For every $p\!\in\! U_{v}$, there exists an $X_i$-compatible holomorphic chart $(U,\varphi)$ around $u_{v}(p)\!\in\! \cZ$ and a sequence of non-zero complex numbers $(t_{a,v,i})_{a\in \N}$ (converging to zero) such that 
\bEq{Convtoze_e}
\tn{(uniformly)  }\lim_{a\lra \infty} t_{a,v,i}^{-1} \, x_{0}\circ \varphi \circ u_{a,v'}|_{K} = x_{0}\circ (\nd_{\cN_i}\varphi (\ze_{v,i}|_{K}))
\eEq
on any compact subset $K \!\subset\! U_{v}$.
\eDf

\noindent
Note that, since $f_a$ Gromov converges to the stable map $f$, we have 
\bEq{Convtoze_e2}
\lim_{a\lra \infty}  x_{i}\circ \varphi \circ u_{a,v'}|_{K} = x_{i}\circ \varphi\circ u_{v}|_{K}\qquad \forall i=0,\ldots,n.
\eEq
uniformly on compact sets.
\begin{remark}
\noindent
For two sequences of non-zero complex numbers $(t_a)_{a\in \N}$ and $(t'_a)_{a\in \N}$, we write 
\bEq{ESeq_e}
(t_a)_{a\in \N}\sim (t'_a)_{a\in \N}\qquad \tn{if}\qquad \lim_{a\lra \infty} t_a/t'_a\!=\!1.
\eEq 
The right-hand side of (\ref{ESeq_e}) defines an equivalence relation on the set of such sequences and we denote the equivalence class of a sequence $(t_a)_{a\in \N}$ by  $[(t_a)_{a\in \N}]$. For  an equivalence class $[(t_a)_{a\in \N}]$ and $t\!\in\!\C^*$, the equation 
$$
t [(t_a)_{a\in \N}]:=[(tt_a)_{a\in \N}]
$$
is well-defined and defines an action of $\C^*$ on the set of equivalence classes. Moreover, the operation of point-wise multiplication/divison  between such sequences 
$$
(t_a)_{a\in \N}\cdot  (t'_a)_{a\in \N} = (t_at'_a)_{a\in \N}
$$
descends to a well-defined multiplication/division operation between the equivalence classes. Condition~(\ref{Convtoze_e}) and the equivalence class of the rescaling sequence $[(t_{a,v,i})_{a\in \N}]$ are independent of the choice of the $X_i$-compatible local chart $(\varphi, U)$; if the limit (\ref{Convtoze_e}) holds in one chart, it will hold in every other chart for the same  $ \ze_{v,i}$. It is also clear from (\ref{Convtoze_e}) that if $(t_{a,v,i})_{a\in \N}$ is a sequence of rescaling parameters associated to $\ze_{v,i}$ and  $(t'_{a,v,i})_{a\in \N}$ is a rescaling sequence associated to $c\ze_{v,i}$, for any $c\!\in\!\C^*$, then
$$
c[(t'_{a,v,i})_{a\in \N}]=[(t_{a,v,i})_{a\in \N}].
$$
\end{remark}

\noindent
The relation between the sets of rescaling parameters $\{t_{a,v,i}\}$ in Definition~\ref{STNS1_dfn}  and $\{\la_a\}$ in (\ref{OSeq_e}) plays an important role in the rest of this chapter.

\bLm{MS_lm}
After passing to a subsequence, every sequence (\ref{OSeq_e}) log Gromov converges and the limit is unique up to equivalence. More specifically,
given (\ref{OSeq_e}), after passing to a subsequence, the limiting holomorphic map $f$ admits meromorphic sections $\{\ze_{v,i}\}_{v\in \V, i\in I_{v}}$ as in Definition~\ref{STNS1_dfn}; furthermore, 
\bEn
\item these meromorphic sections are unique up to multiplication by a constant in~$\C^*$;
\item $\ze_{v,i}$ has no pole/zero in $\Si_{v}\!-\!q_{v}$, 
\item $\ze_{v,i}$ has a zero/pole of order $s_{\uvec{e},i}$ at $q_{\uvec{e}}$, for all $\uvec{e}\in \uvec{\E}_{v}$, $i\!\in\! I_{v}$;
\item for each $\uvec{e}\in \uvec{\E}$, the vector $s_{\uvec{e}}=(s_{\uvec{e},i})_{i \in I_{e}}$ defined as in (\ref{seDf_e}) belongs to $\Z^{I_{e}}_\bullet$;
\item for each $v\in \V$, the product
$$
\bigotimes_{i\in I_{v}} \ze_{v,i} \otimes \bigotimes _{i \in \cI-I_{v}} u_{v}^*\xi_{I_{v},i}
$$ 
is a constant section.
\eEn
\eLm
\bPf
Except items (4) and (5), the rest directly follow from applying \cite[Prp. 3.10 and Lmm.~3.13]{FRelative} to the SNC divisor $\cZ_0\subset \cZ$. 
Also, with respect to the decomposition $\E\cong\E'\cup \E^{\tn{cut}}$, by \cite[Lmm.~3.13]{FRelative}, $s_{\uvec{e}}=0$ for all $e\!\in\! \E'\!\subset\! \E$. Therefore, in order to prove (4), we can restrict to edges in $\E^{\tn{cut}}$. We prove (4) and (5) by following and adjusting the details of \cite[Lmm.~3.13]{FRelative}. \\

\noindent
First, let us recall the setup used in \cite[Sec. 3.4]{FRelative} that we will also use in the rest of this section. For $a$ sufficiently large, the domain $\Si'_a\cong \Si'$ of $f_a$ is obtained from the nodal domain $\Si$ of $f$ in the following way. There exist
\bIt
\item a sequence of complex structures $\mfj_a=(\mfj_{a,v})_{v\in \V}$ on the nodal domain $\Si\!=\!(\Si_{v})_{v\in\V}$ of $f$,
\item a sequence of local $\mfj_{a,v}$-holomorphic coordinates $z_{a,\uvec{e}}\colon \De_{\uvec{e}}\!\lra\! \C$ around $q_{\uvec{e}}\!\in\!\Si_{v}$, for all $v\!\in\!\V$ and $\uvec{e}\!\in\!\uvec{\E}^{\tn{cut}}_{v}$, and
\item  a sequence of non-zero complex numbers $(\ve_{a,e})_{e\in \E^{\tn{cut}}}$ converging to zero,
\eIt
such that 
\bEn
\item\label{glue-domain_it} $(\Si_a,\vec{z_a})$ is isomorphic  to the smoothing of $(\Si,\vec{z}, \mfj_a=(\mfj_{a,v})_{v\in \V})$ defined by
\bEq{GluignRelation_e}
z_{a,\uvec{e}}z_{a,\scz\ucev{e}}=\ve_{a,e} \qquad \forall~e\!\in\!\E^{\tn{cut}},
\eEq
\item the sequence $(\mfj_{a,v})_{a\in \N}$ $C^\infty$-converges to $\mfj_{v}$ for all $v\!\in\!\V$,
\item\label{limitze_e}  the sequence  $(z_{a,\uvec{e}})_{a\in \N}$ $C^\infty$-converges to $z_{\uvec{e}}$, where $z_{\uvec{e}}\colon \De_{\uvec{e}}\!\lra\! \C$ is some fixed local $\mfj_{v}$-holomorphic coordinate  around $q_{\uvec{e}}\!\in\!\Si_{v}$, for all $v\!\in\!\V$ and $\uvec{e}\!\in\!\uvec{\E}^{\tn{cut}}_{v}$.
\eEn

\noindent
With notation as above, for each $e\in\E^{\tn{cut}}$, the union
$$
A_{e}\defeq\! \De_{\uvec{e}} \cup \De_{\scz\ucev{e}}\!=\!
\{(z_{\uvec{e}},z_{\scz{\ucev{e}}})\in \De_{\uvec{e}} \times \De_{\scz\ucev{e}}\colon z_{\uvec{e}} z_{\scz\ucev{e}}=0\}
$$
is a neighborhood of $q_{e}$ in $\Si$.
We orient each circle $\partial\De_{\uvec{e}}$ in the direction of the counter-clock wise rotation in $\De_{\uvec{e}}\!\subset\!\C$.
By (\ref{GluignRelation_e}), the neck region 
$$
A_{a,e}=
\{(z_{a,\uvec{e}},z_{a,\scz{\ucev{e}}})\in \De_{\uvec{e}} \times \De_{\scz\ucev{e}}\colon z_{a,\uvec{e}} z_{a,\scz\ucev{e}}=\ve_{a,e}\}
$$
in $\Si_a$ is a cylinder with two (oppositely oriented) boundary circles
\bEq{partialAe_e}
\partial A_{a,\uvec{e}}\cong\partial\De_{\uvec{e}}\qquad \tn{and}\quad 
\partial A_{a,\scz{\ucev{e}}}\cong\partial\De_{\scz\ucev{e}}.
\eEq
For sufficiently large $a$, $s_{\uvec{e},i}$ is equal to the 
the winding number of $u_a|_{\partial A_{a,\uvec{e}}}$ around the divisor $X_i$; see the proof of \cite[Lmm.~3.13]{FRelative}. \\
\noindent

\noindent
If $I_{e}\!=\!\{i_1,\ldots,i_k\}$, there exists a sufficiently small neighborhood $U$ around $u(q_{e})\!\in\!\cZ$ with coordinates $(x_1,\ldots,x_{n+1})$ such that 
$$
U\cap X_{i_j} = (x_{j}\equiv 0)\qquad \forall j\in \{1,\ldots,k\}
$$
and the projection map $\pi\!\colon\! \cZ\!\lra\! \C$ has the form
\bEq{Productpi_e}
(x_1,\ldots,x_{n+1})\lra \prod_{j=1}^{k} x_{j}.
\eEq
By (\ref{Convtoze_e}) and (\ref{Convtoze_e2}), for sufficiently large $a$,
$$
s_{\uvec{e},i_j}= \tn{winding number around $X_{i_j}$ of } x_{j}\circ u_{a}|_{\partial A_{a,\uvec{e}}}.
$$
Therefore, by (\ref{Productpi_e}) and since $u_a$ has image in $\cZ_{\la_a}$, we have $\sum_{j=1}^k s_{\uvec{e},i_j}=0$\,; i.e. $s_{\uvec{e}}\in \Z^{I_{e}}_\bullet$.\\

\noindent 
Proof of (5) is similar. Since the sections are holomorphic, it is enough to prove (5) on a sufficiently small open set around any point in $\Si_{v}$. Fix $p\!\in\! U_{v}$ and a local coordinate $z$ on a sufficiently small compact disk $K$ around it. If $I_{v}=\{i_1,\ldots,i_k\}$, there exists a sufficiently small neighborhood $U$ around $u(p)\in (X_{I_{v}}\!-\!\partial X_{I_{v}})\!\subset\! \cZ$ with coordinates $(x_1,\ldots,x_{n+1})$ such that 
$$
U\cap X_{i_j} = (x_{j}\equiv 0)\qquad \forall j\in \{1,\ldots,k\}.
$$
and the projection map $\pi\!\colon\! \cZ\!\lra\! \C$ has the form
\bEq{Productpi_e2}
(x_1,\ldots,x_{n+1})\lra \prod_{j=1}^{k} x_{i_j}.
\eEq
On $\cZ$, the product $\otimes_{i \in \cI}\xi_i$ is a section of the trivial line bundle 
\bEq{Trivialization}
\cO_\cZ(\cZ_0)\cong \cZ\times \C.
\eEq
Recall from (\ref{xiChoice_e}) that, when a smoothing $\cZ$ is given, we choose $\xi_i$ such that the projection of $\otimes_{i \in \cI}\xi_i$ to the $\C$-factor in (\ref{Trivialization}) is equal to $\pi$. Therefore, by (\ref{Productpi_e2}), we can choose the local coordinates so that $\xi_{i_j}|_{U}\!=\! x_{j}$ for all $j\!\in\!\{1,\ldots,k\}$. 
By (\ref{Convtoze_e}) and the assumption above, for large $a$,
$$
\la_a =  \bigotimes _{i \in \cI} u_{a,v'}|_{K}^*\xi_{i}=  \bigotimes _{i_j \in I_{v}} x_{j}\circ \varphi\circ u_{a,v'}|_{K}\otimes   \bigotimes _{i \in \cI-I_{v}} u_{a,v'}|_{K}^*\xi_{i}\approx  \bigotimes _{i_j \in I_{v}} t_{a,v,i_j} \ze_{v,i_j} \otimes   \bigotimes _{i \in \cI-I_{v}} u_{a,v'}|_{K}^*\xi_{i}.
$$
Since 
$$
\lim_{a\lra \infty} u_{a,v'}|_{K}^*\xi_{i}= u_{v}|_{K}^* \xi_i \qquad \forall~i\!\in\!\cI-I_{v},
$$
we conclude that 
$$
\lim_{a\lra \infty} \frac{\la_a}{\prod_{i\in I_{v}}t_{a,v,i}} = c_{v}
$$
for some non-zero constant (section) $c_{v}\in \C^*$ and  thus
$$
\bigotimes_{i\in I_{v}} \ze_{v,i} \otimes \bigotimes _{j \in \cI-I_{v}} u_{v}^*\xi_{I_{v},j}=c_{v}.
$$ 
\ePf

\bCr{Log-Limit_Cr}
After passing to a subsequence, every sequence (\ref{OSeq_e}) has a unique limit (\ref{LogLimit_e11}) which belongs to $\cM^{\tn{plog}}_{g,k}(\cZ_0,A)_{\Gamma}$. 
\eCr

\bPf
 On each $\Si_{v}$, we rescale one of $\ze_{v,i}$ (and thus the corresponding sequence $\{t_{a,v,i}\}$) such that 
$$
\bigotimes_{i\in I_{v}} \ze_{v,i} \otimes \bigotimes _{j \in \cI-I_{v}} u_{v}^*\xi_{I_{v},j}=1\qquad \forall v\in \V.
$$ 
Then, by Lemma~\ref{MS_lm}, $f_{\log}$ satisfies all the properties of Definition~\ref{LogTuple_dfn}.
Note that we will then have 
\bEq{NewL_e2}
\lim_{a\lra \infty} \frac{\la_a}{\prod_{i\in I_{v}}t_{a,v,i}} =1\,.
\eEq
We will use (\ref{NewL_e2}) in the proof of the main result below.
\ePf

\bRm{ComapreToRel_rmk}
Since $\cZ_0\subset \cZ$ is an SNC divisor, and $A\cdot \cZ_0=0$, let $\ov\cM_{g,k}^{\log}(\cZ,\cZ_0,A)$ denote the relative (log) moduli space defined in \cite{FRelative} with trivial contact data (with $\cZ_0$) at the marked points. By \cite[Prp.~3.14]{FRelative}, we already know that the unique limit $f_{\log}$ in (\ref{LogLimit_e11}) belongs to $\ov\cM_{g,k}^{\log}(\cZ,\cZ_0,A)$. The linear map $\vr$ in (\ref{vrStar_e}) is the same as the linear map $\vr$ in \cite[(2.28)]{FRelative} but it has a different domain and target. In the following diagram, the first row is (\ref{vrStar_e}), the second row is \cite[(2.28)]{FRelative}, and the vertical maps are the natural inclusion maps.
$$ 
\xymatrix{
\Z^\E\oplus \big(\bigoplus_{v\in \V} \Z^{I_v}\big)^*  \ar[d]\ar^{\vr}[rr]&&
\bigoplus_{e\in \E} \Z_\bullet^{I_e} \ar[d]\\
\Z^\E\oplus \bigoplus_{v\in \V} \Z^{I_v} \ar^{\vr}[rr]&&
\bigoplus_{e\in \E} \Z^{I_e}\,.
}
$$
\eRm

\noindent
The following proposition shows that (\ref{LogLimit_e11}) actually belongs to $\cM^{\log}_{g,k}(\cZ_0,A)_{\Gamma}$. Since $s_{\uvec{e}}\in \Z^{I_e}_{\bullet}$, by (\ref{Constantc_e}), the kernel of the second row is the same as the kernel of the first row; thus, Condition (C1) of Definition~\ref{LogMapVar_dfn} is the same as Condition~(1) of \cite[Dfn.~2.8]{FRelative}. However, the cokernels in each row and thus the groups $\mc{G}$ are different. In order to distinguish the notation, let us denote the group associated to $\Gamma$ in \cite[Dfn.~2.8.(2)]{FRelative} by $\mc{G}^{\tn{rel}}$. The commutative diagram above induces a homomorphism $\mc{G}\lra\mc{G}^{\tn{rel}}$, but this homomorphism does not need to be injective or surjective.

\bPr{Log_prp11} 
Suppose (\ref{OSeq_e}) is a sequence of stable maps in $\ov\cM_{g,k}(\cZ^*,A)$ that log-Gromov converges to (\ref{LogLimit_e11}) in the sense of Definition~\ref{STNS1_dfn}. Then (\ref{LogLimit_e11}) represents an element of $\ov\cM_{g,k}^{\log}(\cZ_0,A)$.
\ePr

\bPf
By Corollary~\ref{Log-Limit_Cr}, we already know that $f_{\log}$ in (\ref{LogLimit_e11}) is a pre-log map.
By Remark~\ref{ComapreToRel_rmk} above, we also know that $f_{\log}$ satisfies Condition (C1) of Definition~\ref{LogMapVar_dfn}. It just remain to show that $f_{\log}$ satisfies Condition (C2) of Definition~\ref{LogMapVar_dfn} as well.\\

\noindent
The proof uses the relation between the following parameters:
\bIt
\item the parameters $\{\la_a\}_{a\in \N}$ in (\ref{OSeq_e});
\item the local holomorphic coordinates $z_{a,\uvec{e}}\colon \De_{\uvec{e}}\!\lra\! \C$ and $z_{\uvec{e}}\colon \De_{\uvec{e}}\!\lra\! \C$ around the nodal points (see the proof of Lemma~\ref{MS_lm});
\item the local coordinates gluing parameters $\{\ve_{a,e}\}_{a\in \N, e\in \E^{\tn{cut}}}$ in (\ref{GluignRelation_e}),
\item the rescaling parameters $\{t_{a,v,i}\}_{a\in \N, v\in \V, i\in I_v}$ in (\ref{Convtoze_e});
\item and, the leading order terms 
$$
0\!\neq\!\eta_{\uvec{e},i}\!\in \!\cN_i|_{u(q_{e})}
$$  
of $f_{\log}$ on $\De_{\uvec{e}}$ defined before (\ref{etae_e}) with respect to $z_{\uvec{e}}$.
\eIt

\noindent
By \cite[Prp.~3.15]{FRelative}, for every oriented edge $\uvec{e}\in\uvec{\E}^{\tn{cut}}$ that goes from $v_1$ to $v_2$, and $i\in I_{e}$, 
\bEn
\item if $i\in I_{v_1}$ and $i\notin I_{v_2}$, we have 
\bEq{ratio_e1Local}  
\lim_{a\lra \infty} t_{a,v_1,i}\,\ve_{a,e}^{s_{\uvec{e},i}}= 
\frac{\eta_{\uvec{e},i}}{\eta_{\scz\ucev{e},i}}\,;
\eEq
\item if $i\!\in\! I_{v_1}\cap I_{v_2}$, we have 
\bEq{ratio_e2Local}  
\lim_{a\lra \infty}\frac{t_{a,v_1,i}\,\ve_{a,e}^{s_{\uvec{e},i}}}{t_{a,v_2,i}} = 
\frac{\eta_{\uvec{e},i}}{\eta_{\scz\ucev{e},i}}\,.
\eEq
\eEn
Additionally, by (\ref{NewL_e2}), we have
$$
\lim_{a\lra \infty} \frac{\la_a}{\prod_{i\in I_{v}}t_{a,v,i}} =1.
$$

\bLm{VertexOrder_Lm}
We can choose the coordinates $\{z_{\uvec{e}}\}_{\uvec{e}\in \uvec{\E}^{\tn{cut}}}$ and $\{z_{a,\uvec{e}}\}_{\uvec{e}\in \uvec{\E}^{\tn{cut}},a\in \N}$ satisfying (\ref{GluignRelation_e}) and item~\ref{limitze_e} after that, and rescalings of $\{\ze_{v,i}\}_{v\in \V, i\in I_v}$ and $(t_{a,v,i})_{a\in \N, v\in \V, i\in I_v}$ such that 
\begin{gather}
\label{ratio-equal_e2}
t_{a,v_1,i}~\ve_{a,e}^{s_{\uvec{e},i}} = t_{a,v_2,i} \qquad \forall~i\!\in\!I_{v_1}\!\cap\!I_{v_2},~a\!>\!\!>\!1,\\
\label{ratio-equal_e1} 
t_{a,v_1,i}~\ve_{a,e}^{s_{\uvec{e},i}} = 1 \qquad \forall~i\!\in\!I_{v_1}\!-\!I_{v_2},~a\!>\!\!>\!1,\\
\label{product_e} 
\prod_{i\in I_{v_1}}t_{a,v_1,i} = \prod_{i\in I_{v_2}}t_{a,v_2,i} \qquad \forall~v_1,v_2\in\V,~a\!>\!\!>\!1.
\end{gather}
\eLm

\bPf
Throughout the proof we assume that the domains $\Si_a$ are smooth; i.e. $\Gamma'$ is a one vertex graph, $\V'=\{v'\}$, and thus $\E=\E^{\tn{cut}}$. The argument in general reduces to this case by focusing on each component of $\Si_a$; see the adjustments at end of the proof of \cite[Prp.~3.14]{FRelative}. We modify a given choice of 
$$
\{z_{\uvec{e}}\}_{\uvec{e}\in \uvec{\E}^{\tn{cut}}},~~\{z_{a,\uvec{e}}\}_{\uvec{e}\in \uvec{\E}^{\tn{cut}},a\in \N},~~\{\ze_{v,i}\}_{v\in \V, i\in I_v},~~(t_{a,v,i})_{a\in \N, v\in \V, i\in I_v},
$$
to another set satisfying  (\ref{ratio-equal_e2})-(\ref{product_e}).\\

\noindent
First, it follows from (\ref{ratio_e1Local}), (\ref{ratio_e2Local}), and (\ref{NewL_e2}), that 
$$
\prod_{i\in I_{e}} \frac{\eta_{\uvec{e},i}}{\eta_{\scz\ucev{e},i}}=1
\qquad \forall~e\!\in\!\E;
$$
i.e. 
$$
\Big(\frac{\eta_{\uvec{e},i}}{\eta_{\scz\ucev{e},i}}\Big)_{i\in I_e}\in (\C^*)^{I_e}_\bullet.
$$
Fix an orientation $O$ on $\E$, and choose some branch 
$$
\eta=\bigoplus_{\uvec{e}\in O} \eta_{e}\!\in\! \bigoplus_{e\in \E} \C^{I_{e}}_{\bullet},\qquad \eta_{e}=\big(-\log (\eta_{\scz\ucev{e},i}/ \eta_{\uvec{e},i})\big)_{i\in I_{e}}\!\in\! \C^{I_{e}}_\bullet\quad \forall~\uvec{e}\!\in \!O,
$$
of the multi-valued function $\log$. By (\ref{NewL_e2}), for each $v\!\in\!\V$ and any $i\!\in\! I_v$, we can replace $\{t_{a,v,i}\}_{a\in \N}$ with another equivalent (in the sense of (\ref{ESeq_e})) sequence  such that 
$$
\la_a=\prod_{i\in I_{v}}t_{a,v,i}\qquad  \forall~a>\!\!>1.
$$
Then we will have 
$$
(t_{a,v,i})_{v\in \V, i\in I_v}\in \Big(\prod_{v\in \V} \C^{I_{v}}\Big)^* \qquad  \forall~a>\!\!>1.
$$
By (\ref{ratio_e1Local}), (\ref{ratio_e2Local}), and definition of $\vr$ in (\ref{vrStar_e}) (via the chosen orientation $O$), we can choose the branches
$$
\xi_{a}=\big((-\log(\ve_{a,e}))_{e\in \E},(-\log(t_{a,v,i}))_{v\in \V, i\in I_{v}}  \big)\in \C^{\E} \oplus \Big(\bigoplus_{v\in \V} \C^{I_{v}}\Big)^*\qquad \forall~a\!\in\!\N
$$
so that 
$$
\lim_{a\lra \infty} \vr_\C(\xi_{a}) =\eta.
$$
By \cite[Lmm.~3.21]{FRelative} applied to $\vr_\C$, there exists a sequence 
$$
(\xi'_a)_{a\in \N}\!\subset\!\C^{\E} \oplus \Big(\bigoplus_{v\in \V} \C^{I_{v}}\Big)^*
$$ 
such that $\vr_\C(\xi_a-\xi'_a)\!=\!0$ for all $a\!\in\!\N$ and the limit $\lim_{a\lra \infty}\xi'_a\!=\!\xi'$ exists.
Taking the exponential of $\xi'$ and $\xi'_a$, we find elements
$$
\big((\al_{e})_{e\in \E}, (\al_{v,i})_{v\in\V,i\in I_{v}}\big),\big((\al_{a,e})_{e\in \E}, (\al_{a,v,i})_{v\in\V,i\in I_{v}}\big)_{a\in \N}\in (\C^*)^{\E}\times \Big(\prod_{v\in \V} (\C^*)^{I_{v}}\Big)^*
$$
such that 
$$
\lim_{a\lra\infty} \big((\al_{a,e})_{e\in \E}, (\al_{a,v,i})_{v\in\V,i\in I_{v}}\big) =\big((\al_{e})_{e\in \E}, (\al_{v,i})_{v\in\V,i\in I_{v}}\big)
$$
and
\begin{gather}
\label{ratio_e2adjusted}
\frac{\big(\al^{-1}_{a,v_1,i} t_{a,v_1,i}\big)~\big(\al_{a,e}^{-1}\ve_{a,e}\big)^{s_{\uvec{e},i}}}
{\big(\al^{-1}_{a,v_2,i}t_{a,v_2,i}\big)} =1 \qquad \forall~i\!\in\!I_{v_1}\!\cap\!I_{v_2},~a\!\in\!\N
,\\
\label{ratio_e1adjusted}  
\big(\al^{-1}_{a,v_1,i} t_{a,v_1,i}\big)~\big(\al_{a,e}^{-1}\ve_{a,e}\big)^{s_{\uvec{e},i}}= 1 \qquad \forall~i\!\in\!I_{v_1}\!-\!I_{v_2},~a\!\in\!\N,\\
\label{ratio_e3adjusted}  
\al_a^{-1}\la_a\defeq  \prod_{i\in I_{v_1}} \al^{-1}_{a,v_1,i} t_{a,v_1,i} = \prod_{i\in I_{v_2}} \al^{-1}_{a,v_2,i} t_{a,v_2,i} \qquad \forall~v_1,v_2\in \V,~a\!\in\!\N.
\end{gather}
By (\ref{ratio_e2adjusted})-(\ref{ratio_e3adjusted}), for $a$ sufficiently large, replacing 
\bIt
\item $\{z_{\uvec{e}}\}_{\uvec{e}\in O}$ with $\{\al_{e}^{-1} z_{\uvec{e}}\}_{\uvec{e}\in O}$, 
\item  $\{z_{a,\uvec{e}}\}_{\uvec{e}\in O}$ with $\{\al^{-1}_{a,e} z_{a,\uvec{e}}\}_{\uvec{e}\in O}$,
\item  $\{\ve_{a,e}\}_{e\in \E}$ with $\{\al^{-1}_{a,e} \ve_{a,e}\}_{e\in \E}$,
\item $(t_{a,v,i})_{v\in \V, i\in I_{v}}$ with $(\al^{-1}_{a,v,i}t_{a,v,i})_{v\in \V, i\in I_{v}}$, and
\item  $(\ze_{v,i})_{v\in \V, i\in I_{v}}$ with $(\al_{v,i}\ze_{v,i})_{v\in \V, i\in I_{v}}$, 
\eIt
we get a new set of representatives satisfying (\ref{ratio-equal_e2}), (\ref{ratio-equal_e1}), and Definition~\ref{PreLogMap_dfn}.(1)  with 
$$
\la(f_{\log})\!\defeq\!\prod_{i\in I_{v}} \al_{v,i}= \lim_{a\lra \infty} \al_a\qquad \tn{for any }v\!\in\!\V. 
$$
\ePf
\noindent
In order to finish the proof of Proposition~\ref{Log_prp11}, by the lemma above, 
the modified set gives us a pre-log map equivalent to $f_{\log}$ that satisfies 
$$
 \frac{\eta_{\uvec{e},i}}{\eta_{\scz\ucev{e},i}}=1\qquad  \forall~e\in \E^{\tn{cut}},~i\!\in\!I_{e}.
$$
By (\ref{RI_rmk}), we conclude that $\tn{ob}_\Gamma(f_{\log})=1\!\in\!\mc{G}$.
\ePf

\noindent
Proposition~\ref{Log_prp11} applies to a sequence of stable maps in $\cZ^*$. For the proof of Theorem~\ref{Compactness_th}, we also need to consider sequences in $\ov\cM_{g,k}^{\log}(X_\eset,A)$ itself. In this case, the ambient smoothing $\cZ$ is not needed and $X_\eset$ can be any d-semistable SNC variety. Suppose  
\bEq{XEsetSeq_e}
f_{a,\log}\defeq\big(u_{a,v'},\ze_{a,v'}=(\ze_{a,v',i})_{i\in I_{v'}},C_{v'}\big)_{v'\in \V'}\qquad a\in \N
\eEq
is a sequence of log maps in $\ov\cM_{g,k}^{\log}(X_\eset,A)$. Since the set of possible decorated graphs for each fixed $(g,k,A)$ is finite, after passing to a subsequence, we may assume that (1) all $f_{a,\log}$ have the same decorated dual graph $\Gamma'\!=\!\Gamma(\V',\E',\L
')$, and (2) the underlying stable maps $f_a\in \ov\cM_{g,k}(X_\eset,A)$ Gromov-converge to a stable map 
$$
f=\big(u_{v},C_{v}\big)_{v\in \V}\in\ov\cM_{g,k}(X_\eset,A)
$$ 
with dual graph $\Gamma=\Gamma(\V,\E,\L)$.
As in Lemma~\ref{MS_lm}, we can find meromorphic sections $\{\ze_{v,i}\}_{v\in \V,i\in I_v}$ such that 
\bEq{LogLimitXeset_e}
f_{\log}\defeq\big(u_{v},\ze_{v}=(\ze_{v,i})_{i\in I_{v}},C_{v}\big)_{v\in \V}\in \cM_{g,k}^{\tn{plog}}(X_\eset,A)_{\Gamma}.
\eEq
Repeating the proof of Proposition~\ref{Log_prp11}, with $\{\la(f_{a,\log})\}_{a\in \N}$ as in Definition~\ref{PreLogMap_dfn}.(1) in place of $\{\la_a\}_{a\in \N}$ in (\ref{OSeq_e}) yields the following. 

\bPr{LogXeset_prp} 
Suppose (\ref{XEsetSeq_e}) is a sequence of stable log maps in $\ov\cM_{g,k}(X_\eset,A)$ that log-Gromov converges to (\ref{LogLimitXeset_e}). The limit (\ref{LogLimitXeset_e}) is unique up to equivalence and represents an element of $\ov\cM_{g,k}^{\log}(X_\eset,A)$.
\ePr

\newtheorem*{proofof-Compactness_th}
{Proof of Theorem~\ref{Compactness_th}}
\begin{proofof-Compactness_th}
Similarly to the classical case, consider the sequential convergence topologies on $\ov\cM_{g,k}^{\log}(\cZ,A)$ or just $\ov\cM_{g,k}^{\log}(X_\eset,A)$ given by Propositions~\ref{Log_prp11}  and~\ref{LogXeset_prp}: a subset $W$ of  the moduli space is closed if every sequence in $W$ has a subsequence with a log-Gromov limit in $W$.  Note that as in \cite[Sec.~5.1]{MS2}, we must show that convergence with respect to the topology defined above is equivalent to log-Gromov convergence. Since the forgetful map $\ov\cM_{g,k}^{\log}(X_\eset,A)\!\to\! \ov\cM_{g,k}(X_\eset,A)$ is finite-to-one and log-Gromov convergence is a lift of the classical Gromov convergence, this property follows from the the corresponding statement for the Gromov convergence topology on $\ov\cM_{g,k}(X_\eset,A)$. In other words, the five axioms\footnote{Even though \cite[Sec. 5.1]{MS2} is about the genus $0$ moduli spaces, the statements used here are valid in all genus.} in \cite[Lmm.~5.6.4]{MS2} lift to sequences in $\ov\cM_{g,k}^{\log}(X_\eset,A)$.\\

\noindent
Suppose $W\!\subset\!\ov\cM_{g,k}(X_\eset,A)~\tn{or}~\ov\cM_{g,k}(\cZ,A) $ is closed and let $W'\!=\!\iota^{-1}(W)$. Let $(f_{a})_{a\in \N}$ be any sequence in $W'$. Its image $(h_a\!=\!\iota(f_{a}))_{a\in \N}$ in $W$ has a subsequence, still denoted by $(h_a)_{a\in \N}$, that Gromov converges to some $h\!\in\!W$. On the other hand, by Proposition~\ref{Log_prp11}  or~\ref{LogXeset_prp}, $(f_{a})_{a\in \N}$ has a subsequence that log-Gromov converges to some $f\!\in\!\ov\cM_{g,k}(X_\eset,A)~\tn{or}~\ov\cM_{g,k}(\cZ,A) $. By definition, we have $\iota(f)=h$, i.e. $f\!\in\!W'$. Therefore, $W'$ is closed. We conclude that $\iota$ is \textbf{continuous}.\\

\noindent
Let $f$ be an arbitrary log map in $\ov\cM_{g,k}^{\log}(X_\eset,A)$ with the decorated dual graph $\Gamma$ and $h\!=\!\iota(f)$ be the underlying stable map in $\ov\cM_{g,k}(X_\eset,A)$.  Let $(U_a)_{a\in \N}$ be a shrinking basis for the (metrizable) topology of  $\ov\cM_{g,k}(X_\eset,A)$ or $\ov\cM_{g,k}(\cZ,A)$  around $h$. Recall that every stable map $h$ admits at most finitely many log lifts $f$, each of which is uniquely specified by the vector decorations on the nodes of its dual graph (i.e. the contact data $s_{\uvec{e}}$ at the nodes $q_{\uvec{e}}$). 
Recall from the proof of  Lemma~\ref{MS_lm} that, for $a$ sufficiently large, by the classical gluing theorem, the domain of every map $h'$ in $U_a$ is obtained from the nodal domain $\Si$ of $h$ by gluing the nodes in a standard way. Furthermore, the image of $h'$ is $C^0$-close to the image of $h$. The dual graph $\Gamma'$ of $h'$ is a contraction of $\Gamma$ in the sense of (\ref{gammaGraph_map}). With these identifications, if $f'$ is a log lift of $h'$ in $U_a$, by its \textbf{decoration type}, we mean 
\bIt
\item the vector decorations $s_{\uvec{e}}$  at its nodes $q_{\uvec{e}}$, together with 
\item the winding\footnote{contact points with $X_\partial$ are among the nodal points and are away from the neck region.} number of $h'$ around $X_i$ along the circles $\partial A_{\uvec{e}}$ (see (\ref{partialAe_e})) on every neck $A_e$ obtained from gluing the node $q_e$ of the domain of $h$; see the proof of Lemma~\ref{MS_lm}. 
\eIt
Thus, we say $f'$ has \textbf{the same decoration type as} $f$ if 
\bEn
\item at every node of the domain of $f'$ the vector decoration $s_{\uvec{e}}$ is the same as the vector decoration at the corresponding node of $f$, 
\item on every neck $A_e$ the winding number of $h'$ around $X_i$ along the circle $\partial A_{\uvec{e}}$ is the same as the tangency order $s_{\uvec{e},i}$ for $f$.
\eEn

\noindent 
For $a$ sufficiently large, define $U'_a$ be the set of elements $f'$ in $\ov\cM_{g,k}^{\log}(X_\eset,A)$ or $\ov\cM_{g,k}^{\log}(\cZ,A)$ whose image $h'$ under $\iota$ lies in $U_a$ and $f'$ has the same decoration type as $f$.
By (1) and (2) above, the restriction of $\iota$ to $U'_a$ is one-to-one. We show that $U'_a$ is open. Let $(f_{b})_{b\in \N}$ be a sequence in the complement of $U'_a$ that log-Gromov converges to $f'$. After possibly passing to a subsequence, we can assume that the underlying sequence of stable maps $(h_{b})_{b\in \N}$ lies either in $U_a$ or its complement $U^c_a$. In the latter case, by definition, $f'$ belongs to the complement of $U_a'$. In the former case, the decoration type of $f'$ (with respect to $f$) will be the same as the decoration type of $f_{b}$ which is, by definition, different from the decoration type of $f$. Therefore, $f'$ belongs to the complement of $U_a'$. We conclude that  $U_a'$ is open. Furthermore, it is easy to see that  $(U'_a)_{a\in \N}$ is a shrinking basis for the topology of $\ov\cM_{g,k}^{\log}(X_\eset,A)$ or $\ov\cM^{\log}_{g,k}(\cZ,A)$ at $f$. Therefore, the log-Gromov topology on $\ov\cM^{\log}_{g,k}(X_\eset)$ or $\ov\cM^{\log}_{g,k}(\cZ,A)$ is first-countable.  \\

\noindent
Hausdorffness is  the consequence of the uniqueness of the limit. If $Y$ is a first-countable topological space and has the property that every convergent sequence has a unique limit then $Y$ is Hausdorff.
Finally, compactness of $\ov\cM_{g,k}^{\log}(X_\eset,A)$ is the consequence of the existence of the limit.\qed

\end{proofof-Compactness_th}

%------------------------------------------------------------------------------------------------------
\section{Comments on deformation theory}\label{Deformation_s}
%------------------------------------------------------------------------------------------------------
In this section, we first calculate the expected dimension of each stratum 
$$
\cM_{g,k}^{\log}(X_\eset,A)_\Gamma\!\subset\!\ov\cM_{g,k}^{\log}(X_\eset,A)
$$ 
and thus identify the virtually main components of $\ov\cM_{g,k}^{\log}(X_\eset,A)$. We then describe the deformation-obstruction exact sequence at any log curve. \\

\noindent
First, let us review the notion of logarithmic tangent bundle and set up the notation.
\noindent
Let $X$ be a smooth holomorphic manifold and $D\!\subset\!X$ be a normal crossings divisor. Around every point $p\!\in\!X$ there exists a chart $\varphi\colon U\lra \C^n$ with coordinates $(x_1,\ldots,x_n)$, with $n\!=\!\dim_\C X$, such that 
$$
\varphi(D\cap U)\!\equiv \!( x_1\cdots x_k\!=\!0)\subset \C^n \quad\tn{for~some}~~~0\leq k\leq n.
$$
In such coordinates, the sheaf $\cT X$ of holomorphic sections of the complex tangent bundle $TX$ is generated by 
$$
\partial_{x_1}, \cdots, \partial_{x_n}
$$
and the \textbf{log tangent sheaf} $\cT X(-\log D)$ is the sub-sheaf  generated by 
$$
\partial^{\log}_{x_1}\defeq x_1\partial_{x_1}, \ldots,\partial^{\log}_{x_k}\defeq x_k \partial_{x_k},
~~~\partial_{x_{k+1}}, \ldots, \partial_{x_n}.
$$
It is dual to the sheaf $\Omega^{1}_X(\log D)$ of meromorphic $1$-forms with at most simple poles along $D_i$.
Since $\cT X(-\log D)$ is locally free, it is the sheaf of holomorphic sections of a holomorphic vector bundle $TX(-\log D)$. The inclusion $\cT X(-\log D)\!\subset\!\cT X$ gives rise to a holomorphic homomorphism
$$
\iota\colon T X(-\log D)\lra T X
$$
which is an isomorphism away from $D$.\\

\noindent
Now, suppose $X_\eset\!=\!\bigcup_{i\in \cI} X_i$ is an SNC K\"ahler variety. For each $i\!\in\!\cI$, let $TX_i(-\log \partial X_i)$ denote the logarithmic tangent bundle of the pair $(X_i,\partial X_i=X_{\partial})$ defined above. If $X_\eset$ is d-semistable,  then (it follows from  \cite[Thm~5.9]{ACGHOSS} that) $X_\eset$ admits a natural holomoprhic vector bundle $T^{\log}X_\eset$ such that 
$$
T^{\log}X_\eset|_{X_i}=TX_i(-\log \partial X_i).
$$
In other words, the collection of logarithmic tangent bundles  \{$TX_i(-\log \partial X_i)\}_{i\in \cI}$ naturally glue along the singular locus $X_\eset$ to define a vector bundle over $X_\eset$ that plays the role of tangent bundle for the central fiber. If $\cZ$ is a semistable smoothing of $\cZ_0=X_\eset$ as in (\ref{pi_e}), then there is an exact sequence 
$$
T^{\log}X_\eset\lra T\cZ(-\log X_\eset)|_{X_\eset} \lra \cO_{X_{\eset}}
$$
meaning that the logarithmic normal bundle of $X_\eset$ in $\cZ$ is the trivial line bundle $\cO_{X_{\eset}}$. Furthermore, via the holomorphic homomorphisms
$$
T\cZ(-\log X_\eset)\lra T\cZ \quad \tn{and}\quad T\C(-\log 0)\lra T\C
$$
the derivative map
$$
\nd \pi \colon T\cZ\lra T\C
$$
lifts to a surjective log derivative map 
$$
\nd^{\log} \pi \colon T\cZ(-\log X_\eset)\lra T\C(-\log 0)
$$
whose kernel over $\la\neq 0$ is $T\cZ_\la$ and over $\la=0$ is $T^{\log}\cZ_0$. In this sense, $T^{\log}X_\eset$ can be considered as the smooth limit of $T\cZ_\la$ when $\la$ converges to $0$. In local coordinates $x\!=\!(x_0,\ldots,x_n)$ such that $\pi\colon \cZ\lra \De$ is given by $x \lra z=x_0\cdots x_k$, we have 
\bEq{ndpiLog_e}
\nd^{\log} \pi|_{\cZ_0} \big(h_1\partial^{\log}_{x_0}+\cdots+h_k\partial^{\log}_{x_k}+h_{k+1}\partial_{x_{k+1}}+ \cdots+h_{n} \partial_{x_n}\big)=\Big(\sum_{i=0}^k h_i\Big)\partial^{\log}_{z}|_{z=0}.
\eEq

\noindent
Given $\cM_{g,k}^{\log}(X_\eset,A)_\Gamma$, with notation as Remark~\ref{vrs_rmk}, recall that both $\tn{K}$ and $\tn{K}_\bullet$ are free $\Z$-modules and 
\bEq{CompK_e}
\dim~\tn{K}_{\R}=\tn{rank}~\tn{K}=\dim~\tn{K}_{\bullet,\R}+1=\tn{rank}~\tn{K}_\bullet+1.
\eEq

\bLm{ExpectedGamma_lmm}
For any admissible decorated dual graph $\Gamma$, the expected complex dimension of $\cM_{g,k}^{\log}(X_\eset,A)_\Gamma$ is 
\bEq{ExpectedGamma_e}
c_1^{T^{\log}X_\eset}(A) + (n-3)(1-g) + k - \tn{rank}~\tn{K}_{\bullet}.
\eEq
\eLm

\bPf
For each $i\!\in\!\cI$,  
$$
\partial X_i\defeq \bigcup_{j\in \cI-i} X_{ij}
$$
is an SNC divisor in $X_i$.
For the given $\Gamma$, fix an arbitrary choice of indices $\{i_v \in I_v\}_{v\in\V}$. For each $v\in \V$, let  $\Gamma_v$ be the one-vertex graph $\{v\}$ with the labeling $I_{v}-\{i_v\}\subset \cI-\{i_v\}$. Also, using the identification map $\Z^\cI_\bullet \cong \Z^{\cI-\{i_v\}}$ in (\ref{pibullet_e}), let $\mfs_v$ be the set of contact vectors in $\Z^{\cI-\{i_v\}}$ at the nodal points $\uvec{\E}_v$ together with the trivial contact vectors at the marked points $\vec{z}_v$. 
For every $f=\big(u_v,\ze_v,C_v\big)_{v\in \V}\in \cM_{g,k}^{\tn{plog}}(X_\eset,A)_\Gamma$, by forgetting the $i_v$-th meromorphic section $\ze_{v,i_v}$ in the $v$-th component $(u_v,\ze_v,C_v)$, we obtain element of the relative log space
$$
\cM_{g_v,\mfs_v}^{\log}(X_{i_v}, \partial X_{i_v},A_v)_{\Gamma_v}= \cM_{g_v,\mfs_v}^{\tn{plog}}(X_{i_v}, \partial X_{i_v},A_v)_{\Gamma_v}.
$$
constructed  in \cite{FRelative}. Therefore, $f$ belongs to the fiber product space
\bEq{FP_e}
\times_{v\in \V}~ \cM_{g_v,\mfs_v}^{\log}(X_{i_v}, \partial X_{i_v},A_v)_{\Gamma_v}
\eEq
where the fiber product is over the evaluation maps into $X_{I_e}\times X_{I_e}$ at the pairs of nodal points $(q_{\uvec{e}},q_{\scz\ucev{e}})$ for all $e\!\in\!\E$.
Then, $\tn{ob}_{\Gamma}$ in (\ref{ObBullet_e}) is a map
$$
\tn{ob}_{\Gamma}\colon \cM_{g,k}^{\tn{plog}}(X_\eset,A)_\Gamma\cong \times_{v\in \V}~ \cM_{g_v,\mfs_v}^{\log}(X_{i_v}, \partial X_{i_v},A_v)_{\Gamma_v} \lra \mc{G}(\Gamma)
$$
such that 
$$
\cM_{g,k}^{\log}(X_\eset,A)_\Gamma=\tn{ob}_{\Gamma}^{-1}(1).
$$

\noindent
By \cite[Prp. 4.8]{FDeformation},  the complex expected dimension of each $\cM_{g_v,\mfs_v}^{\log}(X_{i_v}, \partial X_{i_v},A_v)_{\Gamma_v}$ is  
$$
c_1^{TX_{i_v}(-\log \partial X_{i_v})}(A_v)+(n-3)(1-g_v)+k_v+|\uvec{\E}_v|-|I_v-\{i_v\}|, \quad\tn{where}~~k_v=|\vec{z}_v|.
$$
Since 
$$
c_1^{TX_{i_v}(-\log \partial X_{i_v})}(A_v)=c_1^{T^{\log}X_\eset}(A_v)\qquad \forall~v\!\in\!\V,
$$
the expected dimension of the fiber product  (\ref{FP_e}) is equal to 
\bEq{FPdim_e}
\aligned
&\sum_{v\in \V}\big(c_1^{T^{\log}X_\eset}(A_v)+(n-3)(1-g_v)+k_v+|\uvec{\E}_v|-|I_v| +1 \big)-\sum_{e\in \E} (n-|I_e|+1)=\\
&c_1^{T^{\log}X_\eset}(A)+(n-3)(1-g)+k-|\E|-\sum_{v\in \V}(|I_v|-1)+\sum_{e\in \E} (|I_e|-1).
\endaligned
\eEq
By (\ref{vrStar_e}), 
$$
\tn{rank}~\tn{K} -\tn{dim}_\C(\mc{G})=|\E|+1+\sum_{v\in \V}(|I_v|-1)-\sum_{e\in \E} (|I_e|-1).
$$
The identity (\ref{ExpectedGamma_e}) follows from the second equation in (\ref{FPdim_e}), the last equation, and (\ref{CompK_e}).
\ePf

\bDf{Main_dfn}
We say an admissible decorated dual graph $\Gamma$ is a \textsf{main} graph if 
\bEq{mainGraph_e}
\tn{K}_\bullet=0.
\eEq
\eDf

\noindent
The set of such $\Gamma$ is the same as rigid configurations considered in \cite{ACGS}.

\bCr{Main-Ones}
The moduli space $\ov\cM_{g,k}^{\log}(X_\eset,A)$ has the correct expected dimension and 
the virtually main strata of $\ov\cM_{g,k}^{\log}(X_\eset,A)$ correspond to main graphs. 
\eCr

\noindent
Unlike the classical case, often, $\ov\cM_{g,k}^{\log}(X_\eset,A)$  has many virtually main strata, and they contribute differently to the VFC of $\ov\cM_{g,k}(\cZ_\la,A)$. The degeneration formula (\ref{Formula_e}) describes the weights.  \\

\noindent
Over a smooth target $X$, the deformation-obstruction long exact sequence at a stable marked curve 
$f\!=\![u,C\!=\!(\Si,\vec{z})]\in\!\cM_{g,k}(X,A)$ is the sequence 
\bEq{Def-Obs_LES_e}
0\lra ~\tn{aut}(C) \stackrel{\de}{\lra} 
\Def(u)\lra \Def(f) \lra~\Def(C) \stackrel{\de}{\lra} 
\Obs(u)  \lra  ~\Obs(f) \lra ~0,
\eEq
where 
\bEq{AutDef_e}
\aligned
\tn{Aut}(C)\!=\!H^0_{\dbar}(\Si, T\Si(-\log z))&,\qquad \Def(C)\!=\! H^1_{\dbar}(\Si, T\Si(-\log z)),\\
\Def(u)\!=\!H^0_{\dbar}(\Si, u^* TX)&, \qquad \Obs(u)\!=\!H^1_{\dbar}(\Si, u^*TX),
\endaligned
\eEq
and $T\Si(-\log z)$ is the logarithmic tangent bundle associated to the marked-points divisor $z\subset\Si$..
Alternatively, we may replace $T\Si(-\log z)$ and $TX$ with the corresponding sheaves of holomorphic sections $\cT\Si(-\log z)$ and $\cT X$, respectively, and use \u{c}ech cohomology.
A similar description is feasible when $\Si$ is nodal; see below.
Furthermore, if $u$ is an immersion with normal bundle 
$$
N_u\!=\!u^*TX/\nd u (T\Si)
$$ 
and there are no marked points, then
$$
\Def(f)\!=\!H^0_{\dbar}(\Si, N_u) \qquad \Obs(f)\!=\!H^1_{\dbar}(\Si, N_u).
$$
If $\tn{Obs}(f)\!=\!0$, then a small neighborhood $B(f)$ of $f$ in $\cM_{g,k}(X,A)$ is a smooth orbifold of the expected dimension (\ref{exp-dim_e}); see \cite[Sec.~24.1]{Mirror} and \cite[Rmk.~6.2.1]{FF}. In the following, we outline the generalization of this setup to the case of analytical log maps.

\begin{remark}
In (\ref{Def-Obs_LES_e}), if $u$ is not an immersion or there are marked points, then the cokernel sheaf
\bEq{du_e}
\cN_u=\frac{u^*\cT X}{\nd u (\cT\Si(-\log z))}
\eEq
admits a decomposition
$$
\cN_u= \cN_u^{\tn{free}} \oplus \cN^{\tn{tor}}_u
$$ 
into a direct sum of a torsion free sheaf with the associated holomorphic vector bundle $N_u$  and a skyscraper sheaf such that
$$
\Def(f)= H^0_{\dbar}(\Si, N_u)\oplus H^0(\Si, \cN^{\tn{tor}}_u) \quad \tn{and}\quad \Obs(f)=H^1_{\dbar}(\Si, N_u);
$$
see  \cite[p.~284-285]{ST}.
\end{remark}

\bRm{AlgDO_rmk}  
In the algebraic language, the cohomology groups in (\ref{AutDef_e}) are described as 
$$
\aligned
\tn{Aut}(C)\!=\!\tn{Hom}(\Om^1_{\Si}(\log z), \cO_\Si)&,\quad  \Def(C)\!=\! \tn{Ext}^1(\Om^1_{\Si}(\log z), \cO_\Si),\\
\Def(u)\!=\!\tn{Hom}(u^*\Om^1_X, \cO_\Si)&,\quad \Obs(u)\!=\!\tn{Ext}^1(u^*\Om^1_X, \cO_\Si).
\endaligned
$$
\eRm

\noindent

\noindent
\bLm{LogDef_lm}
Associated with any pre-log map $f\equiv \big(f_{v}\!\equiv\! (u_v,\ze_v,C_v)\big)_{v\in \V}$ (with notation as in (\ref{fplogSetUp_e})), there exists a natural holomorphic\footnote{i.e. the restriction of $\nd (u,\ze)$ to each irreducible component  of $\Si$ is a holomorphic homomorphism.} homomorphism derivative map
$$
T^{\log}\Si (-\log z)\stackrel{\nd (u,\ze)}{\xrightarrow{\hspace*{1cm}}} u^* T^{\log} X_{\eset}
$$
that generalizes the derivative map in (\ref{du_e}).
\eLm 

\bPf
The map $\nd (u,\ze)$ is defined in the following way. Suppose $p\!\in\! \Si_v$ is not a marked or nodal point. Then $u_v(p)\in X_{I_{v}}-\partial X_{I_v}$ and $\ze_{v,i}(p)\neq 0,\infty$ for all $i\in I_v$. If $I_v=\{i_1,\ldots,i_k\}$, a neighborhood $V$ of $u_v(p)$ in $X_\eset$ can be identified with a neighborhood $U$ of $0$ in the affine variety
$$
(x_1\ldots x_k =0)  \subset  \C^{n+1}
$$
so that 
$$
(U\cap \{x_{a}=0\})\!\cong\!(V\cap X_{i_a})\qquad \forall~a=1,\ldots,k\,.
$$
 The coordinates $x_{a}$ also give local trivializations
\bEq{LocalCoordTriv_e}
\cN_{X_{I_v-i_a}}X_{I_v}|_{(X_{I_v}\cap V)}\cong (X_{I_v}\cap V)\times \C,\qquad \forall~a=1,\ldots,k\,;
\eEq
see (\ref{dcN-XD_e}). The sections $\xi_{I_v,i}$ in (\ref{ComaptibleXi_e}) (which are unique up to scaling) give (unique up to scaling) trivializations 
\bEq{Rest_e}
\cO_{X_{I_v}}(X_{I_v+i})|_{X_{I_v}\cap V}\cong (X_{I_v}\cap V)\times \C\qquad \forall~i\!\in\!\mc{I}-I_v.
\eEq
Using (\ref{Rest_e}), the given trivialization in (\ref{PofO_e}) thus gives us a trivialization
\bEq{d-semi-tri_e}
\bigotimes_{i\in I_v} \cN_{X_{I_v-i}}X_{I_v}|_{(X_{I_v}\cap V)}\cong (X_{I_v}\cap V)\times \C.
\eEq
We can choose the coordinates $x_1,\ldots,x_k$ so that the product of the trivializations in (\ref{LocalCoordTriv_e}) is equal to the trivialization in (\ref{d-semi-tri_e}). In such coordinates,
\bEq{LogGenCond_e}
T^{\log}X_\eset|_{(V\cap X_{I_v})}=\lrc{ \sum_{a=1}^k h_a \partial^{\log}_{x_a}+ \sum_{a=k+1}^{n+1} h_a\partial_{x_a}\mid  h_1 +\cdots+h_k=0};
\eEq
see (\ref{ndpiLog_e}).
For a local coordinate $w$ on an open set $\De$ around $p$, with $p=(w=0)$, define 
\bEq{Defduze_e1}
\nd (u,\ze)(\partial_w)=\sum_{a=1}^{k}\frac{\frac{\partial \ze_{v,i_a}}{\partial w}}{\ze_{i_a}} \partial^{\log}_{x_a}+ \sum_{a=k+1}^{n+1} \frac{\partial (x_{a}\circ u_v)}{\partial w}\partial_{x_{a}}.
\eEq
Here $\partial \ze_{v,i}/\partial w$ is defined using the local trivialization (\ref{LocalCoordTriv_e}). It is clear from the definition that (\ref{Defduze_e1}) is invariant under constant rescalings of $\ze_{v,i}$. Also, it follows from Definition~\ref{LogTuple_dfn}(a) and (\ref{d-semi-tri_e}) that (\ref{Defduze_e1}) satisfies (\ref{LogGenCond_e}). If $p$ is one of the marked points, then the local generator of $T^{\log}\Si (-\log z)$ is $\partial^{\log}_w= w\partial w$. Therefore, we define 
$$
\nd (u,\ze)(\partial^{\log}_w)=w \nd (u,\ze)(\partial_w).
$$
It is easy to check that 
$$
\nd_p (u,\ze)\colon T_p^{\log}\Si(-\log z) \lra T_{u_v(p)}^{\log}X_\eset
$$
is independent of the choice of local coordinates $w$ and $x_1,\ldots,x_{n+1}$ used in driving (\ref{Defduze_e1}) and $\nd_p (u,\ze)$ is the zero homomorphism at the marked points.\\

\noindent
Now suppose $q_{e}\!\in\! \Si$ is a nodal point obtained from attaching $\Si_v$ at $q_{\scz\uvec{e}}$ to $\Si_{v'}$ at $q_{\scz\ucev{e}}$. Then 
$$
u(q_e)=u_v(q_{\uvec{e}})= u_{v'}(q_{\scz\ucev{e}})\in X_{I_{e}}-\partial X_{I_e}, \qquad \tn{with} \quad I_e=I_v\cup I_{v'}.
$$ 
Fix local coordinates $w_{\uvec{e}}$ on $\De_{\uvec{e}}$ around $q_{\uvec{e}}$ and $w_{\scz\uvec{e}}$  on $\De_{\scz\uvec{e}}$ around $q_{\scz\uvec{e}}$. Then, as in (\ref{LogGenCond_e}), $T^{\log}\Si(-\log z)$ is generated by $\partial^{\log}_{w_{\uvec{e}}}$ around $q_{\uvec{e}}$ and $\partial^{\log}_{w_{\scz\ucev{e}}}$ around $q_{\scz\ucev{e}}$ satisfying 
\bEq{TlogSiNode_e}
\partial^{\log}_{w_{\uvec{e}}}|_{q_e}=- \partial^{\log}_{w_{\scz\ucev{e}}}|_{q_e}.
\eEq
As before, if $I_e=\{i_1,\ldots,i_k\}$, a neighborhood $V$ of $u(q_{e})$ in $X_\eset$ can be identified with a neighborhood $U$ of $0$ in the affine variety
$$
(x_1\ldots x_k =0)  \subset  \C^{n+1}
$$
such that 
$$
(U\cap \{x_{a}=0\})=(V\cap X_{i_a}), \qquad \forall~a=1,\ldots,k.
$$
As before, the coordinates $x_{a}$ also give local trivializations (\ref{LocalCoordTriv_e}) on $X_{I_v}\cap V$ and $X_{I_{v'}}\cap V$ and the sections $\{\xi_{I_v,i}\}_{i\in \mc{I}-I_e}$ and $\{\xi_{I_{v'},i}\}_{i\in \mc{I}-I_e}$ give compatible (over $X_{I_e}$) local trivializations 
\bEq{Rest_e2}
\cO_{X_{I_v}}(X_{I_v+i})|_{X_{I_v}\cap V}\cong (X_{I_v}\cap V)\!\times\! \C\quad \tn{and}\quad 
\cO_{X_{I_{v'}}}(X_{I_{v'}+i})|_{X_{I_{v'}}\cap V}\cong (X_{I_{v'}}\cap V)\!\times\!\C
\quad \forall~i\!\in\!\mc{I}-I_e.
\eEq
Using (\ref{Rest_e2}), the given trivialization in (\ref{PofO_e}) thus gives us compatible trivializations
\bEq{d-semi-tri_e2}
\Big(\bigotimes_{i\in I_V} \cN_{X_{I_v-i}}X_{I_v}\otimes \bigotimes_{i\in I_e-I_V} \cO_{X_{I_v}}(X_{I_v+i})\Big)|_{(X_{I_v}\cap V)}\cong (X_{I_v}\cap V)\times \C.
\eEq
We can choose the coordinates $x_1,\ldots,x_k$ and local trivializations 
\bEq{Rest_e3}
\aligned
&\cO_{X_{I_v}}(X_{I_v+i})|_{X_{I_v}\cap V}\cong (X_{I_v}\cap V)\times \C\qquad \forall~i\!\in\!I_e-I_v,\\
&\cO_{X_{I_{v'}}}(X_{I_{v'}+i})|_{X_{I_{v'}}\cap V}\cong (X_{I_{v'}}\cap V)\times \C\qquad \forall~i\!\in\!I_e-I_{v'}
\endaligned
\eEq
so that 
$$
\ze_{I_v,i_a}=x_{a} \quad \forall~i_a\in I_{e}-I_v\qquad \tn{and}\qquad \ze_{I_v,i_a}=x_{a} \quad \forall~i_a\in I_{e}-I_{v'},
$$
with respect to (\ref{Rest_e3}), and the trivializations 
\bEq{d-semi-tri_e3}
\aligned
\bigotimes_{i\in I_v} \cN_{X_{I_v-i}}X_{I_v}|_{(X_{I_v}\cap V)}\cong (X_{I_v}\cap V)\times \C,\\
\bigotimes_{i\in I_{v'}} \cN_{X_{I_{v'}-i}}X_{I_{v'}}|_{(X_{I_{v'}}\cap V)}\cong (X_{I_{v'}}\cap V)\times \C,
\endaligned
\eEq
obtained from (\ref{d-semi-tri_e2}) and (\ref{Rest_e3}) coincide with the product trivializations given by  (\ref{LocalCoordTriv_e}).
Then, define
\bEq{Defduze_e2}
\nd (u,\ze)(\partial_{w_{\uvec{e}}}^{\log})= 
\sum_{i_a\in I_{v}}w_{\uvec{e}}\frac{\frac{\partial \ze_{v,i_a}}{\partial w_{\uvec{e}}}}{\ze_{i_a}} \partial^{\log}_{x_a}
+\sum_{i_a\in I_e-I_{v}} w_{\uvec{e}}\,\frac{\frac{\partial (x_{a}\circ u_v)}{\partial w_{\uvec{e}}}}{x_{a}\circ u_v} \partial^{\log}_{x_a}
+\sum_{a=k+1}^{n+1} w_{\uvec{e}}\, \frac{\partial (x_{a}\circ u_v)}{\partial w_{\uvec{e}}}\partial_{x_a}
\eEq
and 
\bEq{Defduze_e3}
\nd (u,\ze)(\partial_{w_{\scz\ucev{e}}}^{\log})= 
\sum_{i_a\in I_{v'}}w_{\uvec{e}}\frac{\frac{\partial \ze_{v,i_a}}{\partial w_{\uvec{e}}}}{\ze_{i_a}} \partial^{\log}_{x_a}
+\sum_{i_a\in I_e-I_{v'}} w_{\scz\ucev{e}}\,\frac{\frac{\partial (x_{a}\circ u_v)}{\partial w_{\scz\ucev{e}}}}{x_{a}\circ u_v} \partial^{\log}_{x_a}
+\sum_{a=k+1}^{n+1} w_{\scz\ucev{e}}\, \frac{\partial (x_{a}\circ u_v)}{\partial w_{\scz\ucev{e}}}\partial_{x_a}.
\eEq
Putting $w_{\uvec{e}}$ and $w_{\scz\ucev{e}}$ equal to zero in (\ref{Defduze_e2}) and (\ref{Defduze_e3}), respectively, we get 
$$
\nd_{q_{\uvec{e}}} (u,\ze)(\partial_{w_{\uvec{e}}}^{\log})=\sum_{a=1}^k s_{\uvec{e},i_a}\partial^{\log}_{x_a}
\quad\tn{and}\quad 
\nd_{q_{\scz\ucev{e}}} (u,\ze)(\partial_{w_{\scz\ucev{e}}}^{\log})=\sum_{a=1}^k s_{\scz\ucev{e},i_a}\partial^{\log}_{x_a}.
$$
By Definition~\ref{PreLogMap_dfn}.\ref{MatchingOrders_l} and (\ref{TlogSiNode_e}),  the map 
$$
\nd_{q_e} (u,\ze)\colon T^{\log}_{q_{\uvec{e}}}\Si(-\log z) \lra T^{\log}_{u(q_{\uvec{e}})}X_\eset
$$
is well-defined.
\ePf
\noindent
The claim is that similarly to (\ref{Def-Obs_LES_e}), for every stable log marked curve 
$f\equiv \big[f_{v}\!\equiv\! (u_v,\ze_v,C_v)\big)_{v\in \V}\big]\in\!\ov\cM_{g,k}^{\log}(X_\eset,A)$, the long exact cohomology sequence
\bEq{Def-Obs_LES_e2}
0\lra ~\tn{aut}(C) \stackrel{\de}{\lra} 
\Def(u,\ze)\lra \Def(f) \lra~\Def(C) \stackrel{\de}{\lra} 
\Obs(u,\ze)  \lra  ~\Obs(f) \lra ~0,
\eEq
arising  (in the sense of Remark~\ref{AlgDO_rmk}) from the short exact sequence of sheaves 
\bEq{ShortLog_e}
\cT^{\log}\Si(-\log z) \stackrel{\nd (u,\ze)}{\xrightarrow{\hspace*{1cm}}} u^* \cT^{\log} X_{\eset} \lra \cN_{u,\ze}= \frac{u^* \cT^{\log} X_{\eset}}{\nd (u,\ze) \big(\cT^{\log}\Si(-\log z)\big)}
\eEq 
on $\Si$ is the deformation-obstruction long exact sequence at $f$.  It also follows from (\ref{ShortLog_e}) and (\ref{Def-Obs_LES_e2}) that the expected $\C$-dimension of $\ov\cM_{g,k}(X_\eset,A)$ is
$$
c_1^{T^{\log}X_\eset}(A)+(n-3)(1-g)+k.
$$

\begin{remark}
An element $\xi\!\in\!\Def(u,\ze)$ is a continuous section of $u^* T^{\log} X_\eset$ such that $\xi_v=\xi|_{\Si_v}$ is a holomorphic section of the vector bundle $u_v^*TX_i(-\log \partial X_i)$ for all $v\!\in\! \V$ and any $i\!\in\! I_v$.
While the map $\nd (u,\ze)$ is defined for arbitrary pre-log map,  by the continuity of $\xi$ at the nodes, the deformation space $\Def(u,\ze)$ only consists of those infinitesimal deformations of $(u,\ze)$ that preserve Conditions (C1) and (C2) in Definition~\ref{LogMapVar_dfn}.
\end{remark}

%------------------------------------------------------------------------------------------------------
\section{The degeneration formula}\label{Formula_s}
%------------------------------------------------------------------------------------------------------

In this section, we describe an explicit (degeneration) formula for the contributions of the virtually main components of $\ov\cM^{\tn{log}}_{g,k}(\cZ_0,A)$ in Corollary~\ref{Main-Ones} to its (hypothetical) VFC.\\

\noindent
Assume $X_\eset\!=\!\bigcup_{i\in \cI} X_i$ is a $d$-semistable SNC K\"ahler variety.  
For each $g,k\!\in\!\N$ and $A\!\in\!H_2(X_\eset,\Z)$, the moduli space $\ov\cM_{g,k}^{\log}(X_\eset,A)$ decomposes into a union of (virtually) main components 
\bEq{ModuliDecomposition_e}
\ov\cM_{g,k}^{\log}(X_\eset,A)= \bigcup_{\tn{main}~\Gamma} \ov\cM_{g,k}^{\log}(X_\eset,A)_{\Gamma}.
\eEq
By Definition~\ref{Main_dfn} and (\ref{CompK_e}), if  $\Gamma$  is a main graph,  up to scaling, there is a unique pair of functions $(s\colon\!\V\!\to\!\R^\cI,\la\colon\!\E\!\to\!\R_+)$ satisfying  Definition~\ref{LogMapVar_dfn}.(C1). The condition~(\ref{mainGraph_e})  implies that the image of the dual $\Z$-linear map 
\bEq{DtoTDual_e2}
\T^\vee  \stackrel{\vr_\bullet^\vee}{\xrightarrow{\hspace*{1.5cm}}} \D_\bullet^\vee
\eEq
is  a sub-lattice of finite index. Let 
$$
m(\Gamma)\defeq \abs{\D_\bullet^\vee/\tn{Im}(\vr^\vee_\bullet)}\in\Z_+.
$$
The degeneration formula (\ref{Formula_e}) below claims that the decomposition (\ref{ModuliDecomposition_e}) gives rise to a similar but weighted decomposition formula for VFCs in a semistable degeneration.

\begin{claim}[Prospect Degeneration Formula]\label{Deg-Formula_th}
Let $\pi\colon\cZ\!\lra\!\De$ be a semistable degeneration as in (\ref{pi_e}). Then, 
for every $A\!\in\!H_2(\cZ_0,\Z)$ and $g,k\!\in\!\N$, we have
\bEq{Formula_e}
[\ov\cM_{g,k}(\cZ_\la,A)]^\vfc=   \sum_{\tn{main } \Gamma}  \frac{m(\Gamma)}{|\aut(\Gamma)|}~[\ov\cM_{g,k}^{\log}(X_\eset,A)_\Gamma]^{\vfc}\;,
\eEq
where $|\aut(\Gamma)|$ is the order of the automorphism group of the decorated dual graph $\Gamma$.
\end{claim}

\noindent
The equality (\ref{Formula_e}) should be thought of as an equality of  \v{C}ech cohomology classes in $\ov\cM^{\log}_{g,k}(\cZ,A)$ in the sense of  \cite[Rmk.~8.2.4]{MW}.

\begin{lemma}
In the case of basic degenerations, i.e. $\cI\!=\!\{1,2\}$, this formula coincides with the Jun Li's formula.
\end{lemma} 
\bPf
It is easy to see that the only decorated graphs with $ \tn{K}_\bullet(\Gamma)\!=\!0$ are bipartite graphs with one set of vertices $\V_1\!\subset\! \V$ indexed by $\{1\}$ and the opposite set $\V_2\!\subset \!\V$ indexed by $\{2\}$. Let us orient the edges to go from $\V_1$ to $\V_2$. Write $s_{\uvec{e}}\!=\!(-m_{\uvec{e}},m_{\uvec{e}})\!\neq\! 0 \in \Z^{2}_\bullet$ for all $e\!\in\! \uvec{\E}_{\V_1,\V_2}$. For each $v\in \V$, since $|I_v|=1$, we have $\Z^{I_v}_\bullet=0$.  For each $e \in \E$, we have $I_e=\{1,2\}$. By picking the second coordinate, we can identify $\Z^{I_{e}}_\bullet$ with $\Z$. Therefore, we have
$$
\vr_\bullet\colon \D_\bullet\cong \Z^\E\lra \T\cong \Z^\E, \qquad (\la_e)_{e\in \E}\lra 
(\la_em_e)_{e\in \E }.
$$ 
We conclude that $m(\Gamma)=\prod_{e\in \E}m_e$ as in Jun Li's formula.
\ePf
\noindent
Unlike the basic case above, a main decorated dual graph for an SNC variety with non-trivial 3-fold (and higher) strata may have components mapped into
a stratum $X_I$ with $|I|\!\geq\!2$; see Section~\ref{Cubic_s}.\\

\noindent
A main step in establishing (\ref{Formula_e})  is to prove a gluing theorem for smoothing the nodes of a log map $f$ as in Definition~\ref{LogMapVar_dfn} to get maps in $\cZ_{\la}$ with $\la\neq 0$. The space of gluing parameters for a fixed log map $f$ with the decorated dual graph $\Gamma$ is a sufficiently small neighborhood of the origin in 
\bEq{GluignEquation_e2}
\aligned
\cN_\Gamma=\bigg\{\big((\ve_e)_{e\in \E},(t_{v,i})_{v\in \V, i\in I_v}\big)&\!\in\!\C^\E\times \prod_{v\in \V}\C^{I_v}\colon 
\prod_{i\in I_v} t_{v,i}=\prod_{i\in I_{v'}} t_{v',i}~~~\tn{and}~~~ \ve_e^{s_{\uvec{e},i}}t_{v,i}=t_{v',i}\\
&\forall~v,v'\!\in\!\V,~e\!\in\!\E_{v,v'},i\!\in\!I_e,~ \uvec{e}~~\tn{s.t}~~s_{\uvec{e},i}\geq 0\bigg\}\subset \C^\E\times \prod_{v\in \V}\C^{I_v}.
\endaligned
\eEq 
In (\ref{GluignEquation_e2}), if $i\!\in\!I_e\!-\!I_v$, by $t_{v,i}$ we mean $1$.
The complex numbers $\ve_e$ are the gluing parameters for the nodes of $\Si$ and $t_{v,i}$ are the parameters for pushing $u_v$ out in the direction of $\ze_{v,i}$. The common value $\la\!=\!\prod_{i\in I_v} t_{v,i}$ describes the fiber $\cZ_\la$ that will contain the glued map. In other words, the projection map $\pi\!\colon\!\cN_\Gamma\!\lra\!\C$ induced by $\pi\!\colon\!\cZ\!\lra\!\C$ is the map
$$
 \big((\ve_e)_{e\in \E},(t_{v,i})_{v\in \V, i\in I_v}\big)\lra \la= \prod_{i\in I_v} t_{v,i}.
$$
\noindent
Let 
$$
\vr^\vee\colon \T^\vee\lra \D^\vee
$$
denote the dual of (\ref{vrStar_e}). The image of $\vr^\vee$ is a finite index sub-lattice of $\tn{K}^\vee$. Define 
$$
m_{\tn{red}}(\Gamma)= \abs{\tn{K}^\vee/\tn{Im}(\vr^\vee)}\in\Z_+.
$$
\bLm{GluignSpace_prp}
The space of gluing parameters $\cN_\Gamma$ is a possibly non-irreducible and non-reduced affine toric sub-variety of $\C^\E\times \prod_{v\in \V} \C^{I_v}$  that is isomorphic to $m_{\tn{red}}(\Gamma)$ copies of the irreducible reduced affine toric variety $Y_{\si(\Gamma)}$ (see after (\ref{Constantc_e})), counting with multiplicities.
In particular, if $\Gamma$ is a main graph, then $Y_\si\cong \C$ and $\pi\!\colon\!\cN_\Gamma\!\lra\!\C$ is a map of degree 
$$
m(\Gamma)=m_{\tn{red}}(\Gamma)\cdot m_\si,
$$
where $m_\si$ is the degree of $\pi\colon Y_\si\lra \C$.
\eLm

\noindent
\bPf
Proof of the first claim is identical to the proof of \cite[Prp.~5.7]{FRelative}. We skip it here.
For the second part, if $\Gamma$ is a main graph, from the dual of the commutative diagram (\ref{CD_e}), we get the exact sequence
$$
0\lra \tn{ker}(\vr^\vee) \lra \tn{Ker}(\vr_\bullet^\vee) \stackrel{\de}{\lra}\Z^\vee\cong \Z\lra \frac{\D^\vee}{\tn{Im}(\vr^\vee)} \lra \frac{\D_\bullet^\vee}{\tn{Im}(\vr_\bullet^\vee)} \cong \frac{\Z}{m(\Gamma) \Z}\, .
$$
The image of $\de$ can not be the entire $\Z$. Therefore, it should be zero, and we get the exact sequence
$$
\Z\lra \frac{\D^\vee}{\tn{Im}(\vr^\vee)} \lra  \frac{\Z}{m(\Gamma) \Z}\, .
$$
On one hand, the degree of $\pi\!\colon\!\cN_\Gamma\!\lra\!\C$ is the index of the first inclusion map above which, by the exactness of the sequence, is equal to $m(\Gamma)$.
On the other hand, by the first statement of the lemma, this degree is the product of the multiplicity $m_{\tn{red}}(\Gamma)$ and the degree $m_\si$ of $\pi\colon Y_\si\lra \C$.
\ePf

\bRm{ACG_rmk}
Except for the coefficients $m(\Gamma)$, the degeneration formula (\ref{Formula_e}) coincides with \cite[(1.1.1)]{ACGS}.
The multiplicity $m_\tau$ in \cite[(1.1.1)]{ACGS} is the integer $m_{\sigma}$ in Lemma~\ref{GluignSpace_prp}. In the case of basic degenerations, the relation between \cite[(1.1.1)]{ACGS} and Jun Li's formula is explained in \cite{KLR}. 
\eRm

%------------------------------------------------------------------------------------------------------
\section{Rational curves in a pencil of cubic surfaces}\label{Cubic_s}
In this section, we re-study the example of the degeneration  of degree $3$ rational curves in a pencil of cubic surfaces, originally studied in \cite[Sec.~6]{ACGS}. As our calculations show, we can easily identify the space of such log maps without blowing up the triple intersection.  \\

\noindent
Let $P$ be a homogenous cubic polynomial in $x_0,\ldots,x_3$ and
$$\cZ'=\!\big\{\big(\la,[x_0,x_1,x_2,x_3]\big)\!\in\!\C\!\times\!\P^3\!:
x_1x_2x_3=\la P(x_0,x_1,x_2,x_3)\big\}\subset\C\!\times\!\P^3\,.$$
Let $\pi'\!:\cZ'\!\lra\!\C$ be the projection map to the first factor. 
For a generic $P$ and $\la\!\neq\!0$, 
$\pi'^{-1}(\la)$ is a smooth cubic hypersurface (divisor) in~$\P^3$. 
For $\la\!=\!0$, $\pi'^{-1}(0)$ is the SNC variety 
\begin{gather*}
X_{\eset}'=\{0\}\!\times\!\big(X_1'\!\cup\!X_2'\!\cup\!X_3'\big)\subset \{0\}\!\times\!\P^3
\qquad\hbox{with}\\
X_i'\!\equiv\!(x_i=0)\approx\P^2~~\forall~i\!\in\!\{1,2,3\},  \quad 
X_{ij}'\!\equiv\!X_i'\!\cap\!X_j'\approx\P^1~~\forall~i,j\!\in\!\{1,2,3\},~i\!\neq\!j.
\end{gather*}
However, the total space $\cZ'$ of $\pi'$ is not smooth at the 9 points of 
$$\cZ'_{\sing}\equiv
\{0\}\!\times\!\big(X_{\partial}'\!\cap\!(P\!=\!0)\big)\subset X_{\eset},
\qquad\hbox{where}\quad
X_{\partial}'=X_{12}'\!\cup\!X_{13}'\!\cup\!X_{23}'\subset\P^3\,.$$
A small K\"aher resolution $\cZ$ of $\cZ'$ can be obtained by blowing up each singular point 
on $X_{ij}'$ in either\footnote{Not all choices result in a K\"ahler configuration.} $X_i'$ or~$X_j'$. 
The map~$\pi'$ then induces a projection $\pi\!:\cZ\!\lra\!\De$ and defines 
a semistable degeneration.
Every fiber of~$\pi$ over~$\C^*$ is a smooth cubic surface. 
The central fiber $\pi^{-1}(0)$ is the SNC variety 
$X_{\eset}\!\equiv\!X_1\!\cup\!X_2\!\cup\!X_3$ with 3~smooth components (i.e. $N\!=\!3$), 
each a blowup of~$\P^2$ at some number of points.
If each singular point on $X_{ij}'$ is blown up in $X_i'$ with $i\!<\!j$,
then $\cZ$ is obtained from~$\cZ'$ through two global blowups of $\C\!\times\!\P^3$
and is thus projective; see Figure~\ref{CubicToB3P2_fg}.\\

\begin{figure}
\begin{pspicture}(-5,-2.6)(11,1.8)
\psset{unit=.3cm}
\psline[linewidth=.07](0,-2)(8,-2)
\psline[linewidth=.07](0,-2)(0,6)
\psline[linewidth=.07](0,-2)(-5.34,-7.66)
\pscircle*(4,-2){.2} \pscircle*(3,-2){.2}\pscircle*(5,-2){.2}
\pscircle*(0,1){.2} \pscircle*(0,2){.2}\pscircle*(0,3){.2}
\pscircle*(-2.9,-5){.2} \pscircle*(-3.8,-6){.2}\pscircle*(-1.9,-4){.2}
\rput(4.5,2.5){${X}'_1$}
\rput(-5.5,0){${X}'_2$}
\rput(2,-6){${X}'_3$}
\rput(9.2,-2.1){${X}'_{31}$}
\rput(0.1,6.8){${X}'_{12}$}
\rput(-5.8,-8.1){${X}'_{23}$}
\rput(1.3,-1.2){${X}'_{123}$}

\psline[linewidth=.07](20,-2)(28,-2)
\psline[linewidth=.07](20,-2)(20,6)
\psline[linewidth=.07](20,-2)(14.34,-7.66)
\pscircle*(24,-2){.2} \pscircle*(23,-2){.2}\pscircle*(25,-2){.2}
\psline[linewidth=.07](24,-2)(24,0) \psline[linewidth=.07](23,-2)(23,0)\psline[linewidth=.07](25,-2)(25,0)
\pscircle*(20,1){.2} \pscircle*(20,2){.2}\pscircle*(20,3){.2}
\psline[linewidth=.2](20,1)(18,-3) %Alpha line  
\psline{<-}(18.8,-.8)(16,2)%Pointer to Alpha line  
\rput(15.3,2.5){$\al$}%Label of Alpha line  
\psline[linewidth=.2](20,1)(22,1) \psline[linewidth=.07](20,2)(22,2)\psline[linewidth=.07](20,3)(22,3)
\pscircle*(17,-5){.2} \pscircle*(16,-6){.2}\pscircle*(18,-4){.2}
\psline[linewidth=.07](17,-5)(17,-3) \psline[linewidth=.07](16,-6)(16,-4)\psline[linewidth=.07](18,-4)(18,-2)
\rput(24.5,2.5){$X_1$}
\rput(15.5,0){$X_2$}
\rput(22,-6){$X_3$}
\rput(29.2,-2.1){$X_{31}$}
\rput(20.1,6.8){$X_{12}$}
\rput(13.8,-8.1){$X_{23}$}

\end{pspicture}
\caption{On left, the central fiber of $\cZ'$ with its $9$ singular points. On right, the central fiber of $\cZ$ with its $9$ exceptional curves. }
\label{CubicToB3P2_fg}
\end{figure}

\noindent
For the count of degree $1$ and degree $2$ rational curves in $\cZ_\la$, it can be shown that all the limiting curves are of the classical type (i.e. they do not pass through the triple intersection). For example, the broken curve $\al$ in Figure~\ref{CubicToB3P2_fg} is one of the $27$ degree $1$ rational curves in the limit. 
For each $\la\!\neq\!0$, the moduli space $\ov\cM_{0,2}(\cZ_\la,[3])$ of $2$-marked genus $0$ degree $3$ maps in $\cZ_\la$ is of the (expected) complex dimension $4$. In degree $3$, for generic $\la$, there are $84$ such curves passing through $2$ generic points of $\cZ_\la$ at the marked points. 
In the limiting SNC variety $X_\eset$, assuming that the two point constraints move to $X_1$ and $X_2$, $81$ of these $84$ maps can be identified among the maps that do not intersect $X_{123}$. There is, however, a new type of main graph $\Gamma$ contributing to the degeneration formula  (\ref{Formula_e}) that has no analogue in the Jun Li's formula. We are going to describe this $\Gamma$, identify the space of log maps $\cM^{\log}_{0,2}(X_\eset,[3])_\Gamma$, and calculate the coefficient $m(\Gamma)$.\\

\noindent
Let $\Gamma$ be the graph with the set of vertices $\V\!=\{v_0,v_1,v_2,v_3\}$ and the set of edges $\E\!=\!\{e_1,e_2,e_3\}$ such that $e_i$ connects $v_0$ and $v_i$, for all $i\!=\!1,2,3$. Choose the orientations $\uvec{e}_i$ to end at $v_0$,  for all $i\!=\!1,2,3$, and assume 
$$
\aligned
&I_{v_0}\!=\!\cI\!=\!\{1,2,3\},\quad s_{\uvec{e}_1}=(-2,1,1)\!\in\!\Z^3_\bullet,\quad s_{\uvec{e}_2}=(1,-2,1)\!\in\!\Z^3_\bullet,\quad s_{\uvec{e}_3}=(1,1,-2)\!\in\!\Z^3_\bullet, \\
 &I_{v_i}=\{i\}, \qquad A_{v_i}=[1]\!\in\!H_2(X_i,\Z), \quad  \forall~i\!=\!1,2,3,
 \endaligned
$$
where, for each $i\!=\!1,2,3$, $[1]\!\in\!H_2(X_i,\Z)$ is the pre-image  in $X_i$ of the class of a line in $X'_i$ away from the blow-up points.
The two legs corresponding to the two marked points are attached to $v_1$ and $v_2$; see Figure~\ref{star_fg}. A log curve with this dual graph is made of 
\bIt
\item a line $\ell_i\!=\!\tn{Im}(u_{v_i})$ in $X_i'\cong \P^2$ passing though the point\footnote{The choice of resolution and the exceptional curves are irrelevant in the following calculations.} $X'_{123}$ for each $i\!\in\!\{1,2,3\}$, and
\item a log tuple 
$$
\Big(u_{v_0}, \{\ze_j\}_{j\in \{1,2,3\}}, \Si_{v_0}\cong \P^1, q_{v_0}=\{ q_{\scz\ucev{e}_j}\}_{j\in \{1,2,3\}}\Big)
$$ 
such that 
\bEn 
\item $u_{v_0}$ is the constant map onto $X_{123}$, and
\item each $\ze_{j}$ is a meromorphic section of the trivial bundle 
$$
u_{v_0}^*\cN_{j}\cong \Si_{v_0}\!\times\!\C
$$
with a zero of order $2$ at $q_{\scz\ucev{e}_j}$ and poles of order $1$ at $\{ q_{\scz\ucev{e}_k}\}_{k\in \{1,2,3\}-\{j\}}$.
\eEn
\eIt
The function $s\colon\!\V\!\lra\!\R^3$ given by
$$
 s_{v_1}=(3,0,0),\quad s_{v_2}=(0,3,0), \quad s_{v_3}=(0,0,3),\quad\tn{and}\quad  s_{v_0}=(1,1,1)
$$
satisfies (C1) of Definition~\ref{LogMapVar_dfn} and is the unique such function up to rescaling. Therefore, $\Gamma$ is a main dual graph (i.e. $\tn{K}_\bullet\!=\!0$ or $ \tn{K}\!\cong \!\Z$).
Since the domain and the target of the injective map (\ref{vrStar_e}) are $6$-dimensional and its kernel is $1$-dimensional, we conclude that the obstruction group $\mc{G}$ is $1$-dimensional. In fact, it is isomorphic to $\C^*$. \\

\begin{figure}
\begin{pspicture}(1,-2.5)(11,1.3)
\psset{unit=.3cm}
\psline[linewidth=0.07](20,-2)(28,-2)
\psline[linewidth=0.07](20,-2)(20,6)
\psline[linewidth=0.07](20,-2)(14.34,-7.66)
\psline[linewidth=.07](24,-2)(24,-1) \psline[linewidth=.07](23,-2)(23,-1)\psline[linewidth=.07](25,-2)(25,-1)
\psline[linewidth=.07](20,1)(21,1) \psline[linewidth=.07](20,2)(21,2)\psline[linewidth=.07](20,3)(21,3)
\psline[linewidth=.07](17,-5)(17,-4) \psline[linewidth=.07](16,-6)(16,-5)\psline[linewidth=.07](18,-4)(18,-3)
\psline[linewidth=0.15](20,-2)(22,2)
\psline[linewidth=0.15](20,-2)(16,1)
\psline[linewidth=0.15](20,-2)(22,-4)
\pscircle*(20,-2){.15}
\rput(24.5,2.5){$X_1$}
\rput(15.5,0){$X_2$}
\rput(22,-6){$X_3$}

\pscircle*(40,-2){.3}\rput(39.7,-1){\small{$v_0$}}
\psline(40,-2)(44,2)\pscircle*(44,2){.3}\rput(45,2.5){\small{$v_1$}}\psline(44,2)(44,4)
\psline[linewidth=0.15]{->}(44,2)(42.5,0.5)\rput(47,.5){\scz{$s_{\uvec{e}_1}=(-2,1,1)$}}

\psline(40,-2)(35,-2)\pscircle*(35,-2){.3}\rput(34,-2.2){\small{$v_2$}}\psline(35,-2)(35,-4)
\psline[linewidth=0.15]{->}(35,-2)(37,-2)\rput(33.5,-1){\scz{$s_{\uvec{e}_2}=(1,-2,1)$}}

\psline(40,-2)(40,-7)\pscircle*(40,-7){.3}\rput(40,-8){\small{$v_3$}}
\psline[linewidth=0.15]{->}(40,-7)(40,-5)\rput(44,-6){\scz{$s_{\uvec{e}_3}=(1,1,-2)$}}

\end{pspicture}
\caption{Dual graph $\Gamma$ and the image of a map belonging to  $\ov\cM_{0,2}(X_\eset,[3])_\Gamma$ in $X_\eset$.}
\label{star_fg}
\end{figure}

\noindent
For each set $I$, let $\Z^I/\Z$ denote the quotient by the diagonal subgroup. We have 
$$
\Z^I/\Z\cong (\Z^I_\bullet)^\vee.
$$
The dual map 
$$
\T^\vee = \bigoplus_{i\in \cI} \bigg(\frac{\Z^{I_{e_i}}}{\Z}\bigg) \cong \bigg(\frac{\Z^{3}}{\Z}\bigg)^3 \stackrel{\vr_\bullet^\vee}{\xrightarrow{\hspace*{1.5cm}}} 
\D_\bullet^\vee\cong \big(\Z^\E\big)^\vee \oplus \bigg(\frac{\Z^{I_{v_0}}}{\Z}\bigg)\cong \Z^3\oplus \bigg(\frac{\Z^{3}}{\Z}\bigg)
$$
in (\ref{DtoTDual_e2}) is given by 
$$
\vr^\vee\big([\eta_1],[\eta_2],[\eta_3]\big)=
\bigg( \!(-2\eta_{11}+\eta_{12}+\eta_{13}), (\eta_{21}-2\eta_{22}+\eta_{23}),(\eta_{31}+\eta_{32}-2\eta_{33}), -\big([\eta_1]+[\eta_2]+[\eta_3]\big)\!\bigg),
$$
where $\eta_i\!=\![\eta_{i1},\eta_{i2},\eta_{i3}]\!\in\!\frac{\Z^{3}}{\Z}$, for any~$i\!\in\!\{1,2,3\}$.
It is straightforward to check that 
$$
\tn{Im}(\vr^\vee)=\lrc{(a,b,c,[x,y,z])\!\in\!\Z^3\!\oplus\! \frac{\Z^3}{\Z}\colon a+b+c \equiv x+y+z ~~\tn{mod}~3}.
$$
Therefore, the quotient group $\D_\bullet^\vee/\tn{Im}(\vr_\bullet^\vee)$ is isomorphic to  $\Z_3$ and  is generated by the class of $(1,0,0, [0,0,0])$; i.e. $m(\Gamma)\!=\!3$.

\bRm{Comapre_rmk}
In this example 
$$
m(\Gamma)=3, \qquad m_\si=3, \qquad m_{\tn{red}}(\Gamma)=1.
$$
In the light of Remark~\ref{ACG_rmk}, this explains why our coefficient $m(\Gamma)\!=\!3$ coincides with the one calculated in  \cite[Sec.~6]{ACGS}.
\eRm

\noindent
In the pre-log space $\cM^{\tn{plog}}_{0,2}(X_\eset,[3])_\Gamma$, the three lines $\ell_1,\ell_2,\ell_3$ are allowed to be any line passing through the point $X_{123}$ with some slope in $\C^*$. However, the condition $\tn{ob}_{\Gamma}(f)\!\in\!\mc{G}\cong\C^*$ in Definition~\ref{LogMapVar_dfn}.(C2) puts a restriction on the set of lines $\ell_1,\ell_2,\ell_3$ that give rise to a log curve.\\

\noindent
For each $i\!\in\!\{1,2,3\}$, the line $\ell_i$ is the completion of the image of a map of the form
$$
\C \lra \C^3, \quad z\lra (x_{ij}(z))_{j=1,2,3}\subset \C^3,~~x_{ii}=0,~\tn{and}~x_{ij}(z)=a_{ij}z,~a_{ij}\!\in\!\C^*,\quad \forall j\!\neq\!i.
$$
It follows from definition of $\tn{ob}_{\Gamma}(f)$ in (\ref{PLtoG_e}) that $\tn{ob}_{\Gamma}(f)\!=\!1$  if and only if for any set of $3$ distinct points $q_{\scz{\ucev{e}}_1}, q_{\scz{\ucev{e}}_2}, q_{\scz{\ucev{e}}_3}\!\in\!\Si_{v_0}= \P^1$ and local coordinates $z_{\scz{\ucev{e}}_1}$, $z_{\scz{\ucev{e}}_2}$, $z_{\scz{\ucev{e}}_3}$ around them, respectively, there exists a set of meromorphic sections $(\ze_i)_{i\in \{1,2,3\}}$ of $\P^1\!\times\!\C$ (holomorphic away from $q_{\scz{\ucev{e}}_1}, q_{\scz{\ucev{e}}_2}, q_{\scz{\ucev{e}}_3}$) such that the product $\ze_1\ze_2\ze_3$ is a constant section and
$$
\ze_i(z_{\scz{\ucev{e}}_j})=a_{ji}z_{\scz{\ucev{e}}_j}^{-1} \quad \forall i\!\in\! \{1,2,3\}, j\!\neq\!i.
$$
A straightforward calculation shows that this is possible if and only if 
$$
\frac{a_{12}}{a_{13}} \frac{a_{31}}{a_{32}} \frac{a_{23}}{a_{21}}\!=\!-1,
$$
i.e. the product of the slopes of $\ell_1$, $\ell_2$, $\ell_3$ (in a certain order) is $-1$. In other words, 
$$
\tn{ob}_\Gamma(f)=-\frac{a_{12}}{a_{13}} \frac{a_{31}}{a_{32}} \frac{a_{23}}{a_{21}}\in \C^*.
$$\\

\noindent
In the degeneration formula (\ref{Formula_e}), imposing two generic point constraints in $X_1$ and $X_2$ on the image of the two marked points
fixes $\ell_1$ and $\ell_2$.  Then the slope condition above fixes $\ell_3$. Therefore, since $m(\Gamma)\!=\!3$ and $\aut(\Gamma)\!=\!1$, the contribution of such a star-shaped log map to the GW count of degree $3$ rational curves in a smooth cubic surface passing through two generic points is $3$. Together with the other $81$ classical-type curves, we recover the $2$-point degree $3$ genus $0$ GW invariant of cubic surface which is $84$. \\

\noindent
We finish with some comments on Question $(3)$  in Page~\pageref{Q3}. 
After removing the trivial component $u_{v_0}\colon\!\Si_{v_0}\!\lra\! X_{123}$, the moduli space $\cM^{\log}_{0,2}(X_\eset,[3])_\Gamma$ decomposes into the relative spaces 
$$
\cM_{0,((0,0),(1,1))}(X_1,X_{1;\partial},[1]),\quad \cM_{0,((0,0),(1,1))}(X_2,X_{2;\partial},[1]),\quad\tn{and}\quad \cM_{0,(1,1)}(X_3,X_{3;\partial},[1]).
$$
Thus, one might still hope to be able to get a decomposition formula in a situation like this. However, in higher dimensions and higher degrees, there seems to be no obvious way to get such a decomposition. The following two examples highlight the issue even further.

\bEx{CubicxP1_eg} 
Consider the family $\cY\!=\!\cZ\times \P^1\!\lra\!\C$, where $\cZ$ is as above and $\pi$ is the lift of the projection map $\pi\colon\!\cZ\!\lra\!\C$.
Let $Y_I\!=\!X_I\!\times\!\P^1$, for all $\eset\!\neq\!I\!\subset\!\{1,2,3\}$.
Consider the same dual graph but with $k\!=\!3$ (i.e. with a third marked point on $\Si_{v_3}$),
$$
A_{v_0}=[0,1]\!\in\!H_2(Y_{123},\Z)\cong H_2( \{\tn{point}\}\!\times\!\P^1,\Z)\cong \{0\}\!\times\! \Z,
$$
and 
$$
A_{v_i}=[1,0]\!\in\!H_2(Y_{i},\Z)\cong H_2( X_i\times \P^1,\Z)\cong  H_2( X_i,\Z)\!\times \!\Z \qquad \forall~i\!\in\!\{1,2,3\};
$$
see Figure~\ref{star_fg2}. 
\begin{figure}
\begin{pspicture}(-3,-3)(11,.5)
\psset{unit=.3cm}
\psline(15,-2)(21,-2)
\psline(15,-2)(12,1)
\psline(15,-2)(11,-5.5)
\psline[linewidth=.15](15,-2)(15,-8)
\psline(15,-8)(21,-8)
\psline[linestyle=dashed,dash=1pt ](15,-8)(12,-5)
\psline(15,-8)(11,-11.5)
\rput(17.5,0.5){\small{$Y_1$}}\rput(11,-2){\small{$Y_2$}}\rput(19,-5){\small{$Y_3$}}
\psline[linestyle=dashed,dash=1pt](15,-3.5)(16.5,-2)\psline[linewidth=.15](16.5,-2)(18,-0.5)\rput(21, -0.7){\scz{$\ell_1=\tn{Im}(u_{v_1})$}}
\psline[linestyle=dashed,dash=1pt](15,-5)(12.2,-4.8)\rput(13.7, -4.2){\scz{$\ell_2$}}
\psline[linewidth=.15](15,-6.5)(18,-9.5)\rput(16.3,-6.6){\scz{$\ell_3$}}
\psarc{<-}(10,-2){5}{10}{90}
\rput(8,3){\scz{$\tn{Im}(u_{v_0})$}}

\end{pspicture}
\caption{The image of a map belonging to  $\ov\cM_{0,3}(Y_\eset,[3,1])_\Gamma$ in $Y_\eset$.}
\label{star_fg2}
\end{figure}
The moduli space $\ov\cM_{0,3}(Y_\eset,[3,1])_\Gamma$ is complex $8$ dimensional with the same contributing factor $m(\Gamma)\!=\!3$ to (\ref{Formula_e}). A smooth fiber $\cY_\la$ of $\cY$ is the product of the smooth cubic surface $\cZ_\la$ and~$\P^1$. Let 
\bEq{productGW_e}
\tn{GW}_{0,3}^{\cY_\la}(\tn{pt},\tn{pt},\al\!\times\!\tn{pt})
\eEq 
be the number of bi-degree $[3,1]$ rational curves in $\cY_\la$ with two point constraints and 
$$
\al\!\times\!\tn{pt}\!\in\!H_2(\cY_\la,\Z)= H_2(\cZ_\la\!\times\!\P^1,\Z),
$$ 
where $\al$ is the homology class of the smoothing of the limiting line shown in Figure~\ref{CubicToB3P2_fg}-Right. Since $m(\Gamma)\!=\!3$ as before, and there is a unique $\Gamma$-type log map in $Y_\eset$ with those constraints, we conclude that the contribution of $\Gamma$-type curves to  (\ref{productGW_e}) is again $3$. 
In examples like this, where there is a non-constant map $u_{v}$ in a stratum $X_{I_v}$ with $|I_v|\!\neq\!1$, for any decomposition of $\ov\cM_{g,k}(X_\eset,A)_\Gamma$  into a fiber product of relative spaces, either (1) $u_{v}$ has to be considered in one of the relative moduli spaces (which normally results in relative spaces with $\mfs\!\notin\!\N^N$), or (2) $(\Si_v,u_v)$ should be removed while its non-trivial GW contribution affects the matching conditions of the remaining parts. The first idea is a motivation behind studying punctured log GW invariants \cite{ACGSII}.
\eEx

\vskip.1in
\noindent
{\it The University of Iowa, Department of Mathematics,~mohammad-tehrani@uiowa.edu}

%------------------------------------------------------------------------------------------------------

\end{document}